\documentclass[12 pt] {article}
\usepackage{WideMargins2}

\def\F{{{\mathbb F}}}
\def\Z{{{\mathbb Z}}}

\def\cj#1{{{\overline{#1}}}}

\def\ang#1{{{\langle #1 \rangle}}}

\let\varep=\varepsilon
\let\x=\times

\let\ee = e
\def\ie{{ \it i.e.}}
\def\ang#1{{{\langle #1 \rangle}}}
\def\jacobi#1#2{{{\left(\frac{#1}{#2}\right)}}}
\def\calT{{{\cal T}}}

\def\calA{{{\cal A}}}

\def\calS{{{\cal S}}}

\title{\Large\bfseries  New Wilson-like theorems arising from Dickson polynomials}
\author{Antonia W.~Bluher\\ National Security Agency \\ tbluher@access4less.net}
\date{}
\begin {document}
\fancytitle

\begin{abstract}
Wilson's Theorem states that
the product of all nonzero elements of a finite field $\F_q$ is $-1$.    
In this article, we define some natural subsets $S \subset \F_q^\x$ ($q$ odd) and find formulas
for the product of the elements of $S$, denoted $\prod S$.  
These new formulas are appealing for the simple,
natural description of the sets $S$, and for the simplicity of the product.
An example is
$\prod\{a \in \F_q^\x : \text{$1-a$ and $3+a$ are nonsquares}\} = 2$ if $q\equiv \pm1 \pmod{12}$,
or $-1$ otherwise.
\end{abstract}

\section{Introduction} \label{sec:intro}
John Wilson (1741--1793) stated that if $p$ is prime, then $p$ divides $(p-1)!+1$, {\it i.e.}, the
product of the nonzero elements of the finite field $\F_p$ is $-1$.
Wilson's Theorem is normally stated for fields of prime order, but its proof easily carries over to 
arbitrary finite fields $\F_q$, where $q=p^n$.  One simply pairs each element of $\F_q^\x$ with its multiplicative inverse, noting
that 1 and $-1$ are the only elements that are equal to their own inverse.  Each pair $\{a,a^{-1}\}$ multiplies
to 1, leaving us with $\prod \F_q^\x = 1 \x -1 = -1$.

Oystein Ore \cite{Ore} provides 
the following history of Wilson's Theorem.
\begin{quotation} In the \emph{Meditationes Algebraicae} by Edward Waring, 
published in Cambridge in~1770, one finds [...] several announcements on the
theory of numbers. One of them is the following: For any prime $p$
the quotient
$$ \frac{1 \cdot 2 \cdot \cdots  \cdot (p-1) + 1}p $$
is an integer.

This result Waring ascribes to one of his pupils John Wilson (1741--1793).
Wilson was a senior wrangler at Cambridge and left the field of mathematics 
quite early to study law. Later he became a judge and was knighted. Waring
gives no proof of Wilson's theorem until the third edition of his 
{\it Meditationes}, which appeared in~1782. Wilson probably arrived at the 
result through numerical computations. Among the posthumous papers of 
Leibniz there were later found similar calculations on the remainders of~$n!$, 
and he seems to have made the same conjecture.  The first proof
of the theorem of Wilson was given by J.\ L.\ Lagrange in a treatise that
appeared in~1770.
\end{quotation}

In \cite{Bluher-dickson}, while investigating permutation properties 
of Dickson polynomials,
we discovered the following ``Wilson-like'' theorems: 
\begin{equation} \prod \{ a \in \F_q^\x : \text{$a$ and $4-a$ are nonsquares} \} = 2 \label{oldWilsonEq1}
\end{equation}
and
\begin{equation} \prod \{ a \in \F_q^\x : \text{$-a$ and $4+a$ are nonsquares} \} 
= \jacobi 2q\cdot  2, 
\label{oldWilsonEq2}
\end{equation}
where $q$ is an odd prime power and 
$\jacobi a q$ denotes the Legendre symbol, defined for $a\in\F_q$ by
$$ \jacobi a q = \begin{cases} 1 & \text{if $a$ is a nonzero square in $\F_q$} \\ 
$-1$ & \text{if $a$ is a nonsquare in $\F_q^\x$} \\ 0 & \text{if $a=0$}. 
\end{cases} 
$$
This article finds more formulas of this type, and again Dickson polynomials
are a key ingredient. See \cite{structure} for a short exposition on these topics.

Define the following sets: 
$$ \calS_k^+ = \left\{ a \in \F_q^\x  : \jacobi {a+k} q = 1 \right\},\qquad 
\calS_k^- = \left\{ a \in \F_q^\x  : \jacobi {a+k} q = -1 \right\},$$
$$ \calS_{k,\ell}^{++}  = \calS_{\ell,k}^{++} = \left\{\, a \in \F_q^\x  : \jacobi {a+k} q = 1, \jacobi {a+\ell} q = 1 \,\right\}, $$
$$ \calS_{k,\ell}^{+-}  = \calS_{\ell,k}^{-+} = \left\{\, a \in \F_q^\x  : \jacobi {a+k} q = 1, \jacobi {a+\ell} q = -1 \,\right\},$$
$$ \calS_{k,\ell}^{--}  = \calS_{\ell,k}^{--} = \left\{\, a \in \F_q^\x  : \jacobi {a+k} q = -1, \jacobi {a+\ell} q = -1 \,\right\}, $$
where $k,\ell \in \F_q$ and $k\ne\ell$.
These definitions implicitly depend on $q$, which is always 
assumed to be an odd prime power.
This article finds simple closed formulas for $\prod \calS_k^+$ and 
$\prod \calS_k^-$ 
(Theorem~\ref{thm:Sk}) and $\prod \calS_{k,\ell}^{\pm,\pm}$ (Section~\ref{sec:SklFormulas})
for all $k$ and $\ell$ ($k\ne \ell$).   

Define a related set $\calT_{j,\ell}^{\varepsilon_1,\varepsilon_2}$ for $j,\ell\in\F_q$,
$\varepsilon_1,\varepsilon_2 \in \{+1,-1\}$, $j+\ell\ne 0$, as follows:
\begin{equation} \calT_{j,\ell}^{\varepsilon_1,\varepsilon_2} = \left \{ a \in \F_q^\x :
\jacobi{j-a} q = \varepsilon_1,\ \jacobi {\ell+a} q = \varepsilon_2 \right\}. \label{TjlDef}
\end{equation}
Note that 
\begin{equation} \calT_{j,\ell}^{\varepsilon_1,\varepsilon_2} = \calS_{-j,\ell}^{\varepsilon\varepsilon_1,\varepsilon_2},\qquad
\text{where $\varepsilon=\jacobi{-1}q$}.  \label{ST}
\end{equation}
Our reason for introducing these sets is that 
$\prod \calS_{k,\ell}^{\pm,\pm}$ requires 
separate formulas for the cases $q \equiv 1 \pmod 4$ and $q \equiv 3 \pmod 4$; whereas 
$\prod \calT_{j,\ell}^{\pm,\pm}$ can be expressed with a single formula.
For example, we will show that $\prod \calS_{-4,0}^{--}=2$
if $q \equiv 1 \pmod 4$, or $-1/2$ if $q \equiv 3 \pmod 4$; whereas
$\prod \calT_{4,0}^{--}=2$ for all odd $q$.
Some sample formulas that we will prove are the following:
\begin{equation}\prod \calT_{2,2}^{-+} = (-1)^{\lfloor (q+3)/8\rfloor}\label{T22} \end{equation}
\begin{equation}\prod \calT_{1,3}^{--} =  \begin{cases} 2 & \text{if $q \equiv \pm 1 \pmod{12} $} \\
-1 & \text{otherwise.} \end{cases} \label{T13}
\end{equation}
If $j$ and $4-j$ are nonzero squares in $\F_q$, then we show in Theorem~\ref{thm:main1}({\it iv}) that
\begin{equation} 
\prod \calT_{j,4-j}^{--} = \jacobi{2+\sqrt{j\,}} q \cdot 2 
= \jacobi{4+2\,\sqrt{4-j\,\,}}q  \cdot 2
= \jacobi 2 q \prod \calT_{4-j,j}^{--}.
\label{biquad}
\end{equation}
Here, $2\pm\sqrt j$ are both squares or both nonsquares since their product 
$4-j$ is assumed to be a square, so the Legendre symbol is well defined, and 
similarly for $4\pm2\,\sqrt{4-j\,\,}$.  We show in Theorem~\ref{thm:main2} and Table~\ref{tableTjl2} that
\begin{equation} \text{
If $j$ is a square and $4-j$ is not, then 
$\prod \calT_{j,4-j}^{--}$ is a square root of $j$.}
\label{sqrtj}
\end{equation}
\begin{equation} \text{
If $4-j$ is a square and $j$ is not, then
$\prod \calT_{j,4-j}^{--}$ is a square root of $4-j$.}
\label{sqrtell}
\end{equation}
\begin{equation} \text{
If $j$ and $4-j$ are nonsquares, then
$\prod \calT_{j,4-j}^{+-}$ is a square root of $j/(4-j)$.}
\label{sqrtellj}
\end{equation}
In each case, a prescription is given that determines
which square root to select. Thus, each formula is determined completely, and
not just up to sign.
In this article we find formulas for all $\prod \calT_{j,\ell}^{\varepsilon_1,
\varepsilon_2}$ ($j\ne -\ell$).
See Tables~\ref{tableTjl1} and~\ref{tableTjl2} for a complete list when
$j+\ell=4$. For arbitrary $j,\ell$ with $j+\ell\ne 0$, 
$\prod \calT_{j,\ell}^{\varepsilon_1,\varepsilon_2}$
can be determined with  our rescaling formula, given in Section~\ref{sec:rescaling}.

As hinted by~(\ref{biquad}),
the above results are related to rational reciprocity. 
For example, consider the case $j=\ell=2$. As is well known (and will be proved in Lemma~\ref{lem:simple}),
2 is a square in $\F_q$ iff $q = \pm 1 \pmod 8$. Assuming $\jacobi2q=1$, (\ref{biquad}) implies that
$\prod \calT_{2,2}^{--}=\jacobi{2+\sqrt 2} q \cdot 2$. On the other hand, 
$\prod \calT_{2,2}^{--}=(-1)^{(q-\varepsilon)/8}\cdot 2$ by Theorem~\ref{thm:main1}({\it ii}),
where $\varepsilon = \jacobi{-1}q$.
This equals 2 if $q \equiv \pm1 \pmod {16}$ 
or $-2$ if $q \equiv 8\pm1 \pmod {16}$. The conclusion is that $2+\sqrt 2$ is a square in $\F_q$ 
if $q\equiv \pm1\pmod{16}$
and a nonsquare if $q \equiv 8 \pm1 \pmod {16}$. Although known (see \cite{Lemmermeyer}), 
it is nonetheless interesting that rational reciprocity arises in this context.  
For more examples of how our work relates to rational reciprocity, see
Section~\ref{sec:reciprocity}.

Of independent interest is our Theorem~\ref{thm:correspondence}, which gives an
explicit bijection between 
sets $\{v,1/v,-v,-1/v\}\subset \mu_{2q-2}\cup\mu_{2q+2}$ and $\F_q$.
The bijection is used to neatly characterize
$\tau\in\F_q$ such that $\jacobi \tau q = \varepsilon_1$, 
$\jacobi {\tau+1}q = \varepsilon_2$
for each pair $(\varepsilon_1,\varepsilon_2)$. 
Several applications of this result are given (see Section~\ref{sec:correspondence}), for example 
a simple proof that $\jacobi 2q = 1$ iff $q \equiv \pm 1 \pmod 8$ and
$\jacobi {-3}q=1$ iff $q \equiv 1 \pmod 3$.

\bigskip
\noindent{\it Notation:} Throughout, $q$ denotes a power of an odd prime, $\F_q$ the field with $q$ elements, 
$\F_q^\x$ its group of nonzero elements under multiplication, and
$\cj \F_q$ its algebraic closure. If $(d,q)=1$ then $\mu_d$ denotes the
$d$-th roots of unity in $\cj \F_q$.
If $c$ is a real number, then $\lfloor c \rfloor$
denotes the largest integer that is $\le c$. 
If $S$ is a subset of $\F_q^\x$, 
then $\prod S$ denotes the product of the elements of $S$, and $|S|$ denotes its cardinality.
By convention, $\prod \emptyset = 1$.
If $S$ and $S'$ are disjoint, then $\prod S \x \prod S' = \prod (S\cup S')$.
We assume basic knowledge of finite fields, as can be found in \cite[Chapter 2]{LN}, 
and properties of the Legendre symbol over prime fields $\F_p$, as can be found in
\cite{Niven-Zuckerman}. It is well known that $\F_q^\x$ is cyclic. If $g$ is a generator, then
$g^i$ is a square if and only if $i$ is even, and $g^{(q-1)/2}=-1$. From this one can deduce the well-known
properties of Legendre symbols 
that $\jacobi a q = a^{(q-1)/2}$ and
$\jacobi a q \jacobi b q = \jacobi {ab} q$ for all $a,b\in\F_q$. 
We often identify the symbols $+$ and $-$ with $1$ and $-1$, respectively.
Thus,  for $\varepsilon_1,\varepsilon_2 \in \{+,-\}$ we have
$$\calS_{k,\ell}^{\varepsilon_1,\varepsilon_2} = \left\{ a \in \F_q^\x :
\jacobi{a+k} q = \varepsilon_1, \quad \jacobi{a+\ell} q = \varepsilon_2 \right\},$$
$$\calT_{j,\ell}^{\varepsilon_1,\varepsilon_2} = \left\{ a \in \F_q^\x :
\jacobi{j-a} q = \varepsilon_1, \quad \jacobi{\ell+a} q = \varepsilon_2 \right\},$$
where $j,k,\ell\in\F_q $, $j+\ell\ne0$, and $k-\ell\ne 0$.
We define sets  that are similar to $\calS_{k,\ell}^{\varepsilon_1,\varepsilon_2}$ but may contain $0$:
\begin{equation} \label{calAdef}
\calA_{k,\ell}^{\varepsilon_1,\varepsilon_2} = \left\{ a \in \F_q :
\jacobi{a+k} q = \varepsilon_1, \quad 
\jacobi{a+\ell} q = \varepsilon_2 \right\}.
\end{equation}
Note that
\begin{equation}
\calA_{k,\ell}^{\varepsilon_1,\varepsilon_2} = \begin{cases}
\calS_{k,\ell}^{\varepsilon_1,\varepsilon_2} \sqcup \{0\}
& \text{if $\jacobi k q=\varepsilon_1$ and $\jacobi \ell q = \varepsilon_2$,}\\
\calS_{k,\ell}^{\varepsilon_1,\varepsilon_2} &\text{otherwise,} \end{cases}
\label{SA}
\end{equation}
where $\sqcup$ denotes the disjoint union, that is, the union of sets that are known to be pairwise disjoint.
We sometimes abbreviate ``if and only if'' by ``iff'' or by the symbol $\iff$.
Because it arises often, we set
\begin{equation}
\varepsilon = \jacobi{-1}q = (-1)^{(q-1)/2} \in \Z. \label{epsDef}
\end{equation}
Starting in Section~\ref{subsec:Dickson}, we use the notation
$$\ang u = u + 1/u.$$
If $\lambda,\beta\in \F_q$, $\lambda \ne 0$, and $S \subset \F_q$, we define
$$\lambda S + \beta = \{\, \lambda s + \beta : s \in S \,\}.$$
One sees easily from the definitions that if $\nu = \jacobi{\lambda} q $ then for $\varepsilon_1,\varepsilon_2 \in \{1,-1\}$,
\begin{equation} \lambda \calA_{k,\ell}^{\varepsilon_1,\varepsilon_2}+ \beta = \calA_{\lambda k-\beta,\lambda\ell-\beta}^{\nu \varepsilon_1,\nu\varepsilon_2},
\label{lambdaA}  \end{equation}
\begin{equation}
\calS_{\lambda k}^{\varepsilon_1} = \lambda \calS_{k}^{\varepsilon_1\nu}, \qquad
\calS_{\lambda k,\lambda \ell}^{\varepsilon_1,\varepsilon_2} = \lambda \calS_{k,\ell}^{\nu\varepsilon_1,\nu\varepsilon_2},\qquad 
\calT_{\lambda j,\lambda \ell}^{\varepsilon_1,\varepsilon_2} = \lambda \calT_{j,\ell}^{\nu\varepsilon_1,\nu\varepsilon_2}. \label{lambdaS}
\end{equation}
This leads to ``rescaling formulas'':
\begin{equation} \prod \calS_{\lambda k ,\lambda \ell}^{\varepsilon_1,\varepsilon_2} = \lambda^s \prod \calS_{k,\ell}^{\nu\varepsilon_1,\nu\varepsilon_2},
\qquad
\prod \calT_{\lambda j ,\lambda \ell}^{\varepsilon_1,\varepsilon_2} = \lambda^t \prod \calT_{j,\ell}^{\nu\varepsilon_1,\nu\varepsilon_2},
\label{rescaling}
\end{equation}
where $s$ is the cardinality of $\calS_{k,\ell}^{\nu\varepsilon_1,\nu\varepsilon_2}$ and
$t$ is the cardinality of $\calT_{j,\ell}^{\nu\varepsilon_1,\nu\varepsilon_2}$. Here $s$ and $t$ will be computed in Section~\ref{sec:correspondence},
resulting in our rescaling formula in Section~\ref{sec:rescaling}. Using this, we can (and will) normalize by a 
scalar multiple to make $\ell-k=4$ and $\ell+j=4$.

The author thanks Art Drisko for thoroughly reviewing this article and providing many thoughtful comments.

\section{Formulas for $\prod \calS_{k}^{+}$ and $\prod \calS_{k}^{-}$} \label{sec:Sk}
\bigskip
This section gives simple formulas for $\prod \calS_k^+$ and $\prod \calS_k^-$ for
all $k \in \F_q$. Lemma~\ref{lem:S0plus} is well known.
\begin{lemma} \label{lem:S0plus} $|\calS_0^+|=|\calS_0^-|=(q-1)/2$,
$\prod \calS_{0}^{+}  = (-1)^{(q+1)/2}$, and \ 
$\prod \calS_{0}^{-}  = (-1)^{(q-1)/2}$. 
\end{lemma}
\begin{proof} $|\calS_0^+|=|\calS_0^-|=(q-1)/2$ because there are exactly $(q-1)/2$ nonzero squares
and the same number of nonsquares in $\F_q$.
Let $g$ be a generator of $\F_q^\x$.  Then $\calS_{0}^{+} = \{g^2,g^4,\ldots,g^{2k}\}$
where $k=(q-1)/2$. The product of these is $g^s$, where $s = 2 + 4 + \ldots + 2k
 = k(k+1)$. Since $g^k=g^{(q-1)/2}=-1$, $g^s = (-1)^{(q+1)/2}$.
Since $\calS_0^+\sqcup \calS_0^-=\F_q^\x$, 
$\prod \calS_{0}^{-} = -1/\prod \calS_{0}^{+} = (-1)^{(q-1)/2}$ by Wilson's Theorem.
\end{proof}

\begin{lemma} $|\calS_1^+| = (q-3)/2$ and $\prod \calS_{1}^{+} = (-1)^{(q-1)/2}/2$. 
\end{lemma}
\begin{proof}
Let $R$ be a complete set of representatives for $\F_q^\x / \{\pm1\}$ such that $1\in R$. 
Now $\calS_{1}^{+} = \left\{a\in\F_q^\x : \jacobi{a+1}q=1\right\} = \{ s^2 - 1 : s \in R,\ s\ne 1 \}$, so
$|\calS_1^+|=|R|-1$ and
\begin{eqnarray*} (-1)^{|R|-1} \prod \calS_{1}^{+} &=& \prod\left\{(1-s)(1+s) : s \in R,\ s\ne1 \right\}  \\
&=& \prod\left\{1-t : t\in\F_q^\x,\ t\ne \pm1 \right\} = \prod \F_q^\x/2 = -1/2.
\end{eqnarray*}
Since $|R|=(q-1)/2$, the result follows. 
\end{proof}

\begin{lemma}  $|\calS_1^-|=(q-1)/2$ and $\prod \calS_{1}^{-} = (-1)^{(q-1)/2}\cdot 2$.
\end{lemma}
\begin{proof} For each $a\in \F_q^\x$, $\jacobi{a+1} q$ is 0, 1, or $-1$; furthermore, 
$\jacobi{a+1} q = 0$ iff $a=-1$.  Thus, $\F_q^\x = \calS_1^+ \sqcup \calS_1^- \sqcup \{-1\}$.
This shows that $|\calS_1^-|=(q-1)-1-|\calS_1^+| = q-2-(q-3)/2=(q-1)/2$, and that
$\prod \calS_{1}^{+} \x \prod \calS_{1}^{-}   \x (-1) = \prod \F_q^\x = -1$, by Wilson's Theorem.
Since $\prod \calS_{1}^{+} = (-1)^{(q-1)/2}/2$, we find that $\prod \calS_{1}^{-} = (-1)^{(q-1)/2}\cdot 2$.
\end{proof}

\begin{theorem} \label{thm:Sk} Let $k\in\F_q$. Then
$$ \prod \calS_{k}^{+} = \begin{cases} 
(-1)^{(q-1)/2}/(2k) & \text{if $k$ is a nonzero square in $\F_q$} \\
(-1)^{(q+1)/2}\cdot 2 & \text{if $k$ is a nonsquare in $\F_q$} \\
 (-1)^{(q+1)/2} & \text{if $k=0$;} \end{cases}$$
$$ \prod \calS_{k}^{-} = \begin{cases} 
(-1)^{(q-1)/2}\cdot 2 & \text{if $k$ is a nonzero square in $\F_q$} \\
(-1)^{(q+1)/2}/(2k) & \text{if $k$ is a nonsquare in $\F_q$} \\
 (-1)^{(q-1)/2}& \text{if $k=0$.} \end{cases}$$
\end{theorem}
\begin{proof}
The result when $k=0$ is contained in Lemma~\ref{lem:S0plus}. If $k\ne 0$, then let $\nu=\jacobi k q$.
By (\ref{lambdaS}) applied with $\lambda$ and $k$ replaced by $k$ and 1, $\calS_k^\nu=k\calS_1^+$ and $\calS_k^{-\nu}=k\calS_1^{-}$. 
Since $|\calS_1^+|=(q-3)/2$, $|\calS_1^-|=(q-1)/2$, and $k^{(q-1)/2}=\nu$, we have
\begin{eqnarray*} \prod \calS_k^{-\nu} &=& \prod (k\calS_1^-) = k^{(q-1)/2} \prod \calS_1^- = \nu  (-1)^{(q-1)/2}\cdot 2,\\
\prod \calS_k^\nu &=& \prod (k\calS_1^+) = k^{(q-3)/2} \prod \calS_1^+ = (\nu/k) (-1)^{(q-1)/2}/2.
\end{eqnarray*}
The result follows.
\end{proof}

\section{Relations between $\prod \calS_{k,\ell}^{++}$, $\prod \calS_{k,\ell}^{+-}$, $\prod \calS_{k,\ell}^{-+}$,
and $\prod \calS_{k,\ell}^{--}$.} \label{sec:relations}
\bigskip

In this section, we show that knowing one of 
$\prod \calS_{k,\ell}^{++}$, $\prod \calS_{k,\ell}^{+-}$, $\prod \calS_{k,\ell}^{-+}$, and
$\prod \calS_{k,\ell}^{--}$ allows to easily compute the other three.

\begin{proposition}  \label{prop:relations}
Let $k,\ell$ be distinct elements of $\F_q$. Then
\begin{equation} \label{Sklplus}
\prod \calS_k^+ = \begin{cases} \prod \calS_{k,\ell}^{++} \x \prod \calS_{k,\ell}^{+-} \x (-\ell) & 
\text{if $\jacobi {k-\ell} q = 1$ and $\ell\ne 0$} \\
\prod \calS_{k,\ell}^{++} \x \prod \calS_{k,\ell}^{+-} & \text{if $\jacobi{k-\ell}q = -1$ or $\ell=0$.} \end{cases}
\end{equation}
\begin{equation} \label{Sklminus}
\prod \calS_\ell^- = \begin{cases} \prod \calS_{k,\ell}^{--} \x \prod \calS_{k,\ell}^{+-} & 
\text{if $\jacobi {\ell-k} q = 1$ or $k=0$} \\
\prod \calS_{k,\ell}^{--} \x \prod \calS_{k,\ell}^{+-} \x (-k) & \text{if $\jacobi{\ell-k}{q}=-1$ and $k\ne0$.} 
\end{cases}
\end{equation}
\begin{equation} \label{Sklminus2}
\prod \calS_k^- = \begin{cases} \prod \calS_{k,\ell}^{--} \x \prod \calS_{k,\ell}^{-+} & 
\text{if $\jacobi {k-\ell} q = 1$ or $\ell=0$} \\
\prod \calS_{k,\ell}^{--} \x \prod \calS_{k,\ell}^{-+} \x (-\ell) & \text{if $\jacobi{k-\ell}{q}=-1$ and $\ell\ne0$.} 
\end{cases}
\end{equation}
\end{proposition}
\begin{proof}
$\calS_k^+$  consists of all $a \in \F_q^\x$ such that $a+k$ is a nonzero square in $\F_q$. 
If $a \in \calS_k^+$, then either $a+\ell=0$, $a+\ell$ is a nonzero square, or
$a+\ell$ is a nonsquare in $\F_q$.  The case $a+\ell=0$ can occur iff $-\ell \in \calS_k^+$,
{\it i.e.}, iff $\jacobi{k-\ell}q=1$ and $-\ell \in \F_q^\x$.  Thus,
\begin{equation} \calS_k^+ = \begin{cases} \calS_{k\ell}^{++} \sqcup \calS_{k\ell}^{+-} \sqcup \{-\ell\}, &
\text{if $\jacobi{k-\ell}{q}=1$ and $\ell\ne0$,} \\
\calS_{k\ell}^{++} \sqcup \calS_{k\ell}^{+-} & \text{if $\jacobi{k-\ell}{q}=-1$ or $\ell=0$.}  \end{cases} \label{Skplus}
\end{equation}
This proves the first formula, and the second is proved similarly, using
\begin{equation} \calS_\ell^- = \begin{cases} \calS_{k\ell}^{--} \sqcup \calS_{k\ell}^{+-} \sqcup \{-k\}, &
\text{if $\jacobi{\ell - k}{q}=-1$ and $k\ne0$,} \\
\calS_{k\ell}^{--} \sqcup \calS_{k\ell}^{+-} & \text{if $\jacobi{\ell-k}{q}=1$ or $k=0$.}  \end{cases} \label{Slminus}
\end{equation}
The third formula is obtained by exchanging $k$ and $\ell$ in the second formula.
\end{proof}

Knowing one of $\prod \calS_{k,\ell}^{++}$, $\prod \calS_{k,\ell}^{+-}$, $\prod \calS_{\ell,k}^{+-}$, and
$\prod \calS_{k,\ell}^{--}$ allows to easily compute the other three.
For example, suppose $\prod \calS_{k,\ell}^{++}$ is known.
Since $\prod \calS_k^+$ is known (by Theorem~\ref{thm:Sk}), 
(\ref{Sklplus}) gives a formula for $\prod \calS_{k,\ell}^{+-}$
in terms of known quantities.  
Next, (\ref{Sklminus}) can be used to solve for $\prod \calS_{k,\ell}^{--}$, 
and finally (\ref{Sklminus2}) can be used to solve for $\prod \calS_{k,\ell}^{-+}$.

\section{A correspondence between $\mu_{2q+2}\cup \mu_{2q-2}$ and $\F_q$}
 \label{sec:correspondence}

The main result of this section is a correspondence $\tau=(1/4)(v-1/v)^2$
between elements $v\in \mu_{2q-2}\cup \mu_{2q+2}$ and elements $\tau \in \F_q$,
where $\mu_d$ denotes the set of $d$th roots of unity in $\cj\F_q$.
As will be shown in Theorem~\ref{thm:correspondence}, this correspondence determines
an explicit bijection between
sets $\{v,1/v,-v,-1/v\}\subset \mu_{2q-2}\cup \mu_{2q+2}$ and $\F_q$, and the multiplicative
order of $v$ determines the square classes of $\tau$ and $\tau+1$.
We apply this to easily compute 
the cardinalities of $\calA_{k,\ell}^{\pm,\pm}$, 
$\calS_{k,\ell}^{\pm,\pm}$, and $\calT_{j,\ell}^{\pm,\pm}$. 
Another application is a simple proof of the formulas for 
$\jacobi 2q $ and $\jacobi{-3}q$.
We remark that the cardinalities of $\calA_{k,\ell}^{\pm,\pm}$ 
are well known in the literature and the computation is based on 
Jacobi sums; see \cite{Crypto88} for a brief history. 
Our computation of $|\calA_{k,\ell}^{\pm,\pm}|$ is very different from those in the literature.
The cardinalities of $\calS_{k,\ell}^{\pm,\pm}$ and $\calT_{j,\ell}^{\pm,\pm}$
are needed to prove our rescaling formula and to obtain a
relation between $\prod \calT_{j,\ell}^{\pm,\pm}$
and $\prod \calT_{\ell,j}^{\pm,\pm}$. (See Section~\ref{sec:rescaling}.)

If $d$ is even, then $v \in \mu_d$ implies that $\{v,1/v,-v,-1/v\} \subset 
\mu_d$, and we call this an orbit. The orbit has exactly four elements 
unless $v^4=1$, in which case it has exactly two elements.  Note that if $v'$ is any element of the orbit,
then $\{\pm v',\pm 1/v' \}=\{\pm v, \pm 1/v\}$.  The choice of which one to call $v$ is arbitrary.

\begin{theorem} \label{thm:correspondence}
The orbits $\{v,1/v,-v,-1/v\}$ of $\mu_{2(q-1)} \cup \mu_{2(q+1)}$ correspond one-to-one with $\tau \in \F_q$.
In this correspondence, $\{v,1/v,-v,-1/v\}$ is sent to $\tau=(1/4)(v-1/v)^2$, and $\tau$ is sent to the orbit
$\{\pm\sqrt{\tau+1}\pm\sqrt{\tau}\}$.  Further, $v^4=1 \iff \tau \in \{-1,0\}$.   If
$\tau \not \in \{-1,0\}$ (and consequently, $v^4 \ne 1$), then  for each pair $a,b \in \{\pm1\}\x \{\pm1\}$,  we have
\begin{equation} \text{$\jacobi\tau q=a$ and $\jacobi{\tau+1}q=b \iff v^{q-ab}=b$.} \label{v_correspondence} \end{equation}
\end{theorem}

\begin{proof}  The theorem claims a bijection between orbits $\{\pm v, \pm 1/v\}\subset \mu_{2q-2}\cup \mu_{2q+2}$ 
and elements of $\tau\in\F_q$,  wherein
$\{\pm v, \pm 1/v\}$ corresponds to $\tau=(1/4)(v-1/v)^2$, and $\tau$ corresponds to the orbit
$\{\pm\sqrt{1+\tau}\pm\sqrt\tau\}$.
Several things need to be demonstrated:  
\begin{enumerate}
\item[(a)] If $\tau \in \F_q$, then $\{\pm\sqrt{1+\tau}\pm\sqrt\tau\}$ really is an orbit;
\item[(b)] If $\{\pm v,\pm 1/v\}$ is an orbit, then $\tau=(1/4)(v-1/v)^2$ is in $\F_q$ and is independent of the choice of
$v$ within the orbit;
\item[(c)] The composition $\tau \mapsto \{\pm\sqrt{\tau+1}\pm\sqrt{\tau}\}$ with $\{\pm v,\pm 1/v\} \mapsto (1/4)(v-1/v)^2$
is the identity map on $\F_q$. 
\item[(d)] The composition $\{\pm v,\pm 1/v\} \mapsto (1/4)(v-1/v)^2$ with
$\tau \mapsto \{\pm\sqrt{\tau+1}\pm\sqrt{\tau}\}$ 
is the identity map on the set of orbits.
\end{enumerate}

(a)\ Let $\tau \in \F_q$. For fixed choices of square roots, let
$v=\sqrt{\tau+1}+\sqrt\tau$. Then $1/v=\sqrt{\tau+1}-\sqrt\tau$, because 
$$\left(\sqrt{\tau+1}+\sqrt{\tau}\right)\left(\sqrt{\tau+1}-\sqrt \tau\right)=(\tau+1)-\tau=1.$$
Thus, $\{\pm \sqrt{\tau+1}\pm \sqrt\tau\}=\{\pm v, \pm 1/v\}$.
Since $(\sqrt{a})^q \in \{\sqrt{a},-\sqrt a\}$ when $a\in\F_q$,
we have $v^q \in \{\pm\sqrt{\tau+1}\pm\sqrt{\tau}\}
=\{v,1/v,-v,-1/v\}$. Then $v^{2q} \in \{v^2,1/v^2\}$, whence $v \in \mu_{2q-2}\cup\mu_{2q+2}$.

(b)\ Let $\{\pm v, \pm 1/v\}$ be an orbit.
It is clear that $(1/4)(v-1/v)^2$ is invariant under $v\mapsto -v$ and $v\mapsto 1/v$, so 
$\tau=(1/4)(v-1/v)^2$ is a well-defined function of the orbit 
(independent of which one is called $v$).
Since $v^{2q} \in \{v^2,v^{-2}\}$, we have
$v^q\in \{\pm v, \pm 1/v \}$.  Thus, $\tau=(1/4)(v^q-1/v^q)^2 = \tau^q$,
showing that $\tau \in \F_q$.

(c)\ We show that if $\tau \in\F_q$ and $v = \sqrt{\tau+1}+\sqrt\tau$ then $\tau=(1/4)(v-1/v)^2$. 
Indeed, $v-1/v=(\sqrt{\tau+1}+\sqrt{\tau})-(\sqrt{\tau+1}-\sqrt{\tau}) = 2\sqrt{\tau}$, so $\tau=(1/4)(v-1/v)^2$.

(d)\  Given $v \in \mu_{2(q+1)}\cup \mu_{2(q-1)}$, let $\tau = (1/4)(v-1/v)^2$, and we claim 
$v\in\{\pm\sqrt{\tau+1}\pm\sqrt{\tau}\}$. 
Since $\tau+1=(1/4)(v-1/v)^2+1=(1/4)(v+1/v)^2$, 
we see that $(v-1/v)/2=\sqrt\tau$, $(v+1/v)/2=\sqrt{\tau+1}$
for the right choices of square root, therefore $v=(v+1/v)/2+(v-1/v)/2 = \sqrt{\tau+1}+\sqrt\tau$.

Now we prove $v^4=1\iff \tau \in \{0,-1\}$. The orbits containing a fourth root of unity are $\{1,-1\}$ and $\{i,-i\}$, where
$\pm i$ are the square roots of $-1$. These correspond to $\tau = (1/4)(v-1/v)^2 \in \{0,-1\}$.

Finally, we show the claimed relations between $\jacobi \tau q$, $\jacobi {\tau+1}q$, and $v$. If $\tau$ and $\tau+1$ are squares
then $v \in \F_q$, so $v^{q-1}=1$. If $\tau+1$ is a square and $\tau $ is a nonsquare, then $v^q=\sqrt{\tau+1}-\sqrt{\tau}=1/v$, so
$v^{q+1}=1$. If $\tau+1$ is a nonsquare and $\tau$ is a square, then $v^q=-(\sqrt{\tau+1}-\sqrt\tau)=-1/v$, so $v^{q+1}=-1$.
If $\tau$ and $\tau+1$ are nonsquares, then $v^q=-v$, so $v^{q-1}=-1$.
\end{proof}

Note that we can write $\F_q$ as a disjoint union:
\begin{equation}\F_q = \{0\} \sqcup \{-1\} \sqcup \calA_{0,1}^{++} \sqcup \calA_{0,1}^{+-} \sqcup \calA_{0,1}^{-+} \sqcup \calA_{0,1}^{--}.
\label{disjointUnion}
\end{equation}
In the correspondence, $\tau = 0$ corresponds to the orbit $\{1,-1\}$,
$\tau = -1$ corresponds to the orbit $\{\pm\sqrt{-1}\}$, and $\tau \in \calA_{0,1}^{\varep_1,\varep_2}$ corresponds to all orbits
$\{\pm v,\pm1/v\}$ such that $v^{q-\varep_1\varep_2}=\varep_2$ and $v^4\ne 1$.
In other words,
\begin{equation} \calA_{0,1}^{\varepsilon_1,\varepsilon_2}
= \left\{\, (v-1/v)^2/4 : v^{q-\varepsilon_1\varepsilon_2}=\varepsilon_2,\ 
v^4\ne 1 \,\right\}. \label{v_epsilon}
\end{equation}

\begin{corollary} \label{cor:v2} For all $\varepsilon_1,\varepsilon_2
\in \{\pm1\}$,
\begin{equation} \calA_{-2,2}^{\varepsilon_1,\varepsilon_2} =
\left\{\, v^2 + v^{-2} :
\text{$v^{q-\varepsilon_1\varepsilon_2} = \varepsilon_2$ and $v^4\ne 1$} \,\right\}. \label{v2}
\end{equation}
\end{corollary}
\begin{proof} Consider what happens to the equality (\ref{v_epsilon}) 
when $\tau$ is mapped to $4\tau+2$.
First, $\tau \in \calA_{0,1}^{\varepsilon_1,\varepsilon_2} \iff 4\tau+2 \in
\calA_{-2,2}^{\varepsilon_1,\varepsilon_2}$. Second,
$\tau = (v-1/v)^2/4 \iff 4\tau+2 = v^2+v^{-2}$. The result follows.
\end{proof}

We give several applications of Theorem~\ref{thm:correspondence},  beginning with a simple proof of the formulas for
$\jacobi 2q$ and $\jacobi{-3}q$.   Lemma~\ref{lem:simple} is well known.

\begin{lemma} \label{lem:simple} Let $q = p^n$, where $p$ is an odd prime. 
Then $-1$ is a square in $\F_q$ if and only if $q\equiv 1 \pmod 4$, 
$2$ is a square in
$\F_q$ if and only if $q \equiv \pm1 \pmod 8$, and $-3$ is a nonzero square 
in $\F_q$
if and only if $q \equiv 1 \pmod 3$.  If $a \in \F_p$ then $\jacobi a q = {\jacobi a p}^n$.
If $\varepsilon = (-1)^{(q-1)/2}$, then $(q-\varepsilon)/4$ is an integer and
\begin{equation} \jacobi 2 q = (-1)^{(q-\varepsilon)/4}. \label{jacobi2q} \end{equation}
\end{lemma}

\begin{proof} For the first assertion,
$\jacobi{-1}q=(-1)^{(q-1)/2}$ and $(q-1)/2$ is even iff $q\equiv 1 \pmod 4$. 
For the penultimate assertion, if $a \in \F_p^\x$, then 
$$\jacobi a q = a^{(q-1)/2}=\left(a^{(p-1)/2}\right)^m = {\jacobi a p}^m,$$
where $m=(p^n-1)/(p-1)=1+p+\cdots+p^{n-1}$. 
Since $m\equiv n \pmod 2$, we have $\jacobi a q = \jacobi a p^n$.

Now we prove that $-3$ is a square in $\F_q \iff q \equiv 1 \pmod 3$.
Set $v=\omega$, where $\omega$ is a primitive cube root of unity, and let $\tau\in\F_q$ correspond
to the orbit of~$v$. 
Then $\tau = (1/4)(v-1/v)^2 = (1/4)(\omega^2 + \omega^{-2}-2) = (1/4)(\omega^2+\omega-2)=-3/4$,
and $\tau+1=1/4$.  By Theorem~\ref{thm:correspondence},
$$\text{$\jacobi {-3}q = 1  \iff \tau$ and $\tau+1$ are squares $\iff\omega^{q-1}=1$}$$
and the latter is equivalent to $3|q-1$.

Next, we prove that 2 is a square in $\F_q \iff q \equiv \pm 1 \pmod 8$.
Let $v = \zeta$, a primitive eighth root of unity, and let $\tau\in\F_q$ correspond to the orbit of~$v$.
Then $\tau = (1/4)(\zeta-1/\zeta)^2=(1/4)(i+1/i-2)=-2/4$ and $\tau+1=2/4$, where $i=\zeta^2$. By Theorem~\ref{thm:correspondence},
$$\jacobi {-2}q=\jacobi 2 q = 1 \iff \zeta^{q-1}=1 \iff 8|(q-1),$$
$$\jacobi {-2}q=-1,\ \jacobi 2 q = 1 \iff \zeta^{q+1}=1 \iff 8|(q+1).$$
Thus, $\jacobi 2q = 1$ iff $q\equiv \pm 1 \pmod 8$.

Finally, we prove (\ref{jacobi2q}). Note that $q\equiv \varepsilon \pmod 4$.
Since $(q-\varepsilon)/4$ is even iff $q \equiv \varepsilon \pmod 8$ iff $q\equiv \pm 1 \pmod 8$ iff $\jacobi 2q=1$,
we see that $\jacobi 2q=(-1)^{(q-\varepsilon)/4}$.
\end{proof}

A second application of Theorem~\ref{thm:correspondence} is the following counting result.  Again, the result is known, but the method of
proof is new.

\begin{proposition} \label{prop:AijCard0}
Let $\calA_{k,\ell}^{\pm,\pm}$ be as in~(\ref{calAdef}).    Then
$|\calA_{0,1}^{++}|=\lfloor (q-3)/4 \rfloor$, $|\calA_{0,1}^{+-}| = \lfloor (q+1)/4\rfloor$, 
$|\calA_{0,1}^{-+}|=\lfloor (q-1)/4 \rfloor$, $|\calA_{0,1}^{--}| = \lfloor (q-1)/4\rfloor$. 
\end{proposition}
\begin{proof} 
The theorem shows that $v \mapsto \tau = (1/4)(v-1/v)^2$ gives a four-to-one map from $(\mu_{2(q+1)}  \cup \mu_{2(q-1)} )\setminus \mu_4$
onto $\F_q \setminus \{0,-1\}$, and that 
\begin{equation} \text{$\tau \in \calA_{0,1}^{\varepsilon_1,\varepsilon_2} \iff
v^{q-\varepsilon_1\varepsilon_2}=\varepsilon_2$ and $v^4\ne1$.}  \label{vSet}
\end{equation}
Each set of $v$'s in~(\ref{vSet}) can be described as $S_{a,b}\setminus \mu_4$, where $S_{a,b}=\{v\in\cj\F_q : v^{q+a} = b \}$ and $a,b\in \{1,-1\}$. 
Note that $|S_{a,b}|=q+a$. Also,
the cardinality of $S_{a,b}\setminus \mu_4$ is divisible by~4, since it consists of disjoint orbits of size~4. Let $i=|S_{a,b}\cap \mu_4|$.
Then $0\le i \le 4$, and $|S_{a,b}|-i = q+a-i$ is divisible by~4.   Since $q+a$ is even, $i$ is even. That is, $i \in \{0,2,4\}$.
If $b=1$ then $1\in S_{a,b}\cap \mu_4$, so $i \in \{2,4\}$ and $|S_{a,b}\setminus \mu_4|/4 \in \{(q+a-2)/4,(q+a-4)/4 \}$.  This number must be an integer,
so $|S_{a,b}\setminus \mu_4|/4 = \lfloor (q+a-2)/4 \rfloor$. Thus, 
$|\calA_{0,1}^{++}|=\lfloor (q-3)/4 \rfloor$ and
$|\calA_{0,1}^{-+}| = \lfloor (q-1)/4 \rfloor$. 
If $b=-1$, then $1\not \in S_{a,b}$, so $i \in \{0,2\}$. Then $|S_{a,b}\setminus \mu_4|/4
\in \{(q+a)/4,(q+a-2)/4\}$. Again, this must be an integer, so $|S_{a,b}\setminus \mu_4|/4 = \lfloor (q+a)/4 \rfloor$.
Thus, $|\calA_{0,1}^{+-}| = \lfloor (q+1)/4 \rfloor$ and $|\calA_{0,1}^{--}| = \lfloor (q-1)/4 \rfloor$.
\end{proof}

\begin{corollary} 
If $a$ is a random nonsquare in $\F_q\setminus\{0,-1\}$ (under the uniform distribution), then $a+1$ is a square
with probability exactly 1/2.  
\end{corollary}
\begin{proof}   The value $a$ is randomly selected from $\calA_{0,1}^{-+} \sqcup \calA_{0,1}^{--}$.
The corollary follows by noting that 
$|\calA_{0,1}^{-+}|=|\calA_{0,1}^{--}|$.
\end{proof}

Let $\varepsilon=\jacobi{-1}q$ and $m=(q-\varepsilon)/4 \in \Z$.
In Proposition~\ref{prop:AijCard0}, it is useful to note that $\lfloor (q+1)/4 \rfloor = m$, $\lfloor (q-3)/4 \rfloor
=m-1$, and $\lfloor (q-1)/4 \rfloor = m + (\varepsilon-1)/2$. Thus, 
\begin{equation} |\calA_{0,1}^{++}|=m-1,\quad |\calA_{0,1}^{+-}| = m,\quad
|\calA_{0,1}^{-+}|= |\calA_{0,1}^{--}| = m + (\varepsilon-1)/2. \label{AijCard}
\end{equation}

\begin{proposition} \label{prop:AijCard}   
Let $k,\ell$ be distinct elements of $\F_q$ and $\nu = \jacobi{\ell-k}q$. Then
$$|\calA_{k,\ell}^{\nu,\nu}| = m-1,\quad  |\calA_{k,\ell}^{\nu,-\nu}| = m,\quad 
|\calA_{k,\ell}^{-\nu,\nu}| = |\calA_{k,\ell}^{-\nu,-\nu}|=m + (\varepsilon-1)/2, $$
where $m=(q-\varep)/4$. Also, if $\mu \in \{1,-1\}$ then
\begin{equation} |\calA_{k,\ell}^{\varep\mu,\mu}| = m-\gamma,\qquad \text{where $\gamma = \begin{cases} 1 & \text{if $\mu = \nu$} \\
0 & \text{if $\mu=-\nu$.} \end{cases}$} \label{AijCard2}
\end{equation}
\end{proposition}
\begin{proof} The result is true when $k=0$ and $\ell=1$ by (\ref{AijCard}).
If $k,\ell$ are arbitrary distinct elements of $\F_q$, then
$(\ell-k) \calA_{0,1}^{\nu \varepsilon_1, \nu \varepsilon_2} - k =
\calA_{k,\ell}^{\varepsilon_1,\varepsilon_2}$
by (\ref{lambdaA}).  
In particular, $\calA_{k,\ell}^{\varepsilon_1,\varepsilon_2}$
and $\calA_{0,1}^{\nu \varepsilon_1,\nu \varepsilon_2}$ have the same cardinality.
The result follows.
\end{proof}

\begin{corollary} \label{cor:SklCard}  Let $j,k,\ell\in\F_q$, $k\ne \ell$,
$j+\ell \ne 0$. Then
\begin{eqnarray*} |\calS_{k,\ell}^{\varepsilon_1,\varepsilon_2}| &=& 
\begin{cases} |\calA_{k,\ell}^{\varepsilon_1,\varepsilon_2}|-1 & \text{if
$\jacobi k q = \varepsilon_1$ and $\jacobi \ell q = \varepsilon_2$,} \\
|\calA_{k,\ell}^{\varepsilon_1,\varepsilon_2}| & \text{otherwise;} \end{cases} \\
|\calT_{j,\ell}^{\varepsilon_1,\varepsilon_2}| &=& 
\begin{cases} |\calA_{-j,\ell}^{\varepsilon\varepsilon_1,\varepsilon_2}|-1 & \text{if
$\jacobi j q = \varepsilon_1$ and $\jacobi \ell q = \varepsilon_2$,} \\
|\calA_{-j,\ell}^{\varepsilon\varepsilon_1,\varepsilon_2}| & \text{otherwise.} \end{cases} \\
\end{eqnarray*}
Here $|\calA_{k,\ell}^{\varepsilon_1,\varepsilon_2}|$ and
$|\calA_{-j,\ell}^{\varepsilon\varepsilon_1,\varepsilon_2}|$ can be
determined using Proposition~\ref{prop:AijCard}.
\end{corollary}
\begin{proof} This follows from (\ref{SA}) and~(\ref{ST}).
\end{proof}

\section{Rescaling} \label{sec:rescaling}

In Sections~\ref{sec:SklFormulas} and~\ref{sec:proofs} we will use Dickson polynomials to find formulas for $\prod \calT_{j,\ell}^{\pm,\pm}$
when $j+\ell=4$.  Lemma~\ref{lem:rescaling} allows us to extrapolate these formulas to arbitrary pairs $j'$ and $\ell'$
with $j'+\ell'\ne 0$. Similar comments apply to $\calS_{k,\ell}^{\pm,\pm}$.

\begin{lemma}  (Rescaling Formula) \label{lem:rescaling}
Let $j',\ell'$ be arbitrary elements of $\F_q$ such that $j'+\ell'\ne0$,
and let
$\varepsilon_1,\varepsilon_2 \in \{1,-1\}$. Further, let
$\lambda = (j'+\ell')/4$,
$\nu = \jacobi {j'+\ell'} q$, $j=j'/\lambda$, and $\ell=\ell'/\lambda$.
(Note that $j+\ell=4$). Then
$$\prod \calT_{j',\ell'}^{\varepsilon_1,\varepsilon_2} =
\prod \calS_{-j',\ell'}^{\varepsilon\varepsilon_1,\varepsilon_2} = 
\lambda^{m-\beta-\gamma} \prod \calS_{-j,\ell}^{\nu \varepsilon\varepsilon_1,\nu\varepsilon_2}
=\lambda^{m-\beta-\gamma} \prod \calT_{j,\ell}^{\nu \varepsilon_1,\nu\varepsilon_2},$$
where
$\varepsilon = \jacobi{-1}q$,  $m=(q-\varepsilon)/4$, 
$$\beta = \begin{cases}  1 & \text{ if $\jacobi{j'}q=\varepsilon_1$, 
$\jacobi{\ell'}q = \varepsilon_2$} \\
0 & \text{otherwise;} \end{cases}
$$
$$\gamma = \begin{cases} 1 &
\text{if $\nu=\varep\varep_1 = \varep_2$ or $-\varep=\nu\varep_1=1$}
\\ 0 & \text{otherwise.} 
\end{cases}
$$
\end{lemma}

\begin{proof}  By (\ref{ST}) and (\ref{lambdaS}), 
$$\calT_{j',\ell'}^{\varepsilon_1,\varepsilon_2} =
\calS_{-j',\ell'}^{\varepsilon \varepsilon_1,\varepsilon_2}
= \lambda \calS_{-j,\ell}^{\nu\varep\varep_1,\nu\varep_2}
= \lambda \calT_{j,\ell}^{\nu\varep_1,\nu\varep_2},$$
so the above sets have the same cardinality $t$, and to prove the result
we need only show that $t=m-\beta-\gamma$.  Indeed,
\begin{eqnarray*} t=|\calS_{-j,\ell}^{\nu\varepsilon\varep_1,\nu\varep_2}| &=& 
|\calA_{-j,\ell}^{\nu\varep\varepsilon_1,\nu\varepsilon_2}|-\beta' 
\qquad\text{by Corollary~\ref{cor:SklCard}} \\
&=& m - \gamma' -\beta' \qquad \text{\quad by Proposition~\ref{prop:AijCard},}
\end{eqnarray*}
where 
$$\beta'=\begin{cases} 1 & \text{if $\jacobi{-j}q=\nu\varep\varep_1$ and $\jacobi\ell q=\nu\varep_2$} \\
0 & \text{otherwise,}  \end{cases} \qquad
\gamma'=\begin{cases} 1 & \text{if $\nu\varepsilon\varep_1 = \nu\varep_2=1$} \\
0 & \text{if $\nu \varep\varep_1=1=-\nu\varep_2$} \\
(1-\varep)/2 & \text{if $\nu\varep\varep_1=-1$.} \end{cases} $$
First, $\beta = \beta'$ because
$\jacobi{-j}q =\nu\varep\varep_1$ iff $\jacobi{j'}q=\varep_1$, and
$\jacobi{\ell}q=\nu\varep_2$ iff $\jacobi{\ell'}q=\varep_2$. To show that
$\gamma=\gamma'$, we will show that $\gamma=1$ iff $\gamma'=1$ when $\varep=1$,
and $\gamma=0$ iff $\gamma'=0$ when $\varep = -1$. Suppose $\varep=1$.
Then $\gamma'=1$ iff $\nu\varep_1=\nu\varep_2=1$ iff $\nu=\varep_1=\varep_2$
iff $\gamma=1$. Next, suppose $\varep=-1$. Then $\gamma=1$ iff
$\nu=-\varep_1=\varep_2$ or $\nu=\varep_1$.  Thus, $\gamma=0$ iff
$\nu=-\varep_1=-\varep_2$ iff $-\nu\varep_1=1=-\nu\varep_2$ iff
$\gamma'=0$.
\end{proof}

\begin{proposition} \label{prop:Tellj} Suppose $j,\ell \in \F_q$ and $j+\ell\ne 0$.
Then for $\mu \in \{1,-1\}$,
$$\prod \calT_{\ell,j}^{\mu,\mu} = 
\begin{cases} \mu\jacobi 2q \jacobi{j+\ell}q \prod \calT_{j,\ell}^{\mu,\mu} & \text{if $\jacobi jq=\jacobi\ell q = \mu$} \\
-\mu\jacobi 2q \jacobi{j+\ell}q \prod \calT_{j,\ell}^{\mu,\mu} & \text{otherwise.}
\end{cases} $$
\end{proposition}

\begin{proof} $\calT_{\ell,j}^{\mu,\mu}=
-\calT_{j,\ell}^{\mu,\mu}$ by (\ref{TjlDef}), so $\prod \calT_{\ell,j}^{\mu,\mu} = (-1)^t \prod \calT_{j,\ell}^{\mu,\mu}$, with
$t = |\calT_{j,\ell}^{\mu,\mu}|$. By Corollary~\ref{cor:SklCard}, $t=|\calA_{-j,\ell}^{\varepsilon\mu,\mu}|-\beta$,
where $\beta = 1$ if $\jacobi j q = \jacobi \ell q = \mu$ and $\beta = 0$ otherwise. By~(\ref{AijCard2}),
$|\calA_{-j,\ell}^{\varepsilon\mu,\mu}| = m-\gamma$, where $\gamma=1$ if $\mu\jacobi{j+\ell}q=1$
and $\gamma=0$ if $\mu\jacobi{j+\ell}q=-1$.
Note that $(-1)^\gamma=-\mu\jacobi{j+\ell}q$. Thus, by (\ref{jacobi2q}),
$(-1)^t=(-1)^{m-\gamma-\beta}= \jacobi2q \x -\mu\jacobi{j+\ell}q(-1)^\beta$.  The result follows.
\end{proof}

\section{Legendre symbol identities} \label{sec:legendre}

This section presents
some results involving Legendre symbols. In an earlier draft of this article,
Propositions~\ref{prop:consequence1} and~\ref{prop:consequence2} were proved using Theorem~\ref{thm:correspondence}, 
but then we found a simpler proof.  The results of this section are 
neither obvious nor deep.

\begin{lemma} \label{lem:abc} If $a,b,c\in\F_q$, $ab\ne0$, and $a^2+b^2=c^2$, then 
$\jacobi{c+a}q = \jacobi 2 q \jacobi {c+b} q \ne 0$.
Also, $\jacobi{c+a}q=\jacobi{c-a}q$ and $\jacobi{c+b}q=\jacobi{c-b}q$.
\end{lemma}

\begin{proof} For the first sentence, it suffices to prove that $2(c+a)(c+b)$ is a nonzero square. 
It is nonzero because $2(c+a)(c+b)(c-a)(c-b)=
2(c^2-a^2)(c^2-b^2)=2b^2a^2\ne0$. It is a square because 
$$(a+b+c)^2 = a^2+b^2+c^2 + 2(ab+ac+bc) = 2c^2 + 2(ab+ac+bc) = 2(a+c)(b+c).$$  
The last sentence follows because 
$(c+a)(c-a)=b^2$ and $(c+b)(c-b)=a^2$ are nonzero squares.
\end{proof}

\begin{proposition} \label{prop:consequence1}
If $\tau \in \F_q$ and $\jacobi {\tau (\tau+1)} q = 1$, then \
$\jacobi {2\tau + 1 + 2 \sqrt{\tau (\tau+1)}} q = \jacobi \tau q = \jacobi{\tau+1}q$.
\end{proposition}
\begin{proof} First, $\jacobi{\tau(\tau+1)}q=1$ implies $\jacobi\tau q = \jacobi{\tau+1}q$. For the remaining
equalities, apply Lemma~\ref{lem:abc} with $a=2\sqrt{\tau(\tau+1)}$, $b=1$, and $c=2\tau+1$.
\end{proof}

\begin{proposition} \label{prop:consequence2} 
Suppose that $\tau \in \F_q$ and that $\tau,1+\tau$ are both nonzero squares. Then
\begin{enumerate} 
\item[({\it i})] $1/\tau$ and $1+1/\tau$ are both nonzero squares. 
\item[({\it ii})] $1\pm1/\sqrt{1+\tau}$ are both squares or both nonsquares,
so the Legendre symbol $\jacobi {1\pm1/\sqrt{1+\tau}}q$ is well defined. 
Similarly, the Legendre symbol $\jacobi{1\pm1/\sqrt{1+1/\tau}}q$ 
is well defined. 
\item[({\it iii})] We have
$$\jacobi {1\pm 1/\sqrt{1+1/\tau}} q 
= \jacobi 2 q  \jacobi {1\pm 1/\sqrt{1+\tau} } q.$$
\end{enumerate}
\end{proposition}
\begin{proof} ({\it i}) follows from $1+1/\tau=(1+\tau)/\tau$. For ({\it ii}) and ({\it iii}), 
apply Lemma~\ref{lem:abc} with $a=\pm 1/\sqrt{1+1/\tau}$, $b=\pm1/\sqrt{1+\tau}$, and $c=1$.
Note that $a^2+b^2=c^2$.
\end{proof}

\begin{lemma} If $j,\ell $ are nonzero squares in $\F_q$ and $j+\ell=4$ then
\begin{equation} \label{jell}
\jacobi{2\pm \sqrt {j\,}}q = \jacobi 2 q \jacobi{2\pm \sqrt{\ell\,}}q.
\end{equation}
\end{lemma}

\begin{proof} Apply Lemma~\ref{lem:abc} with $a=\sqrt {j\,}$, $b=\sqrt{\ell\,}$ and $c=2$.
\end{proof}

\section{Formulas for $\prod \calS_{k,\ell}^{\pm, \pm}$ and $\prod \calT_{j,\ell}^{\pm,\pm}$} 
\label{sec:SklFormulas}
\bigskip
In evaluating $\prod \calS_{k,\ell}^{\varepsilon_1,\varepsilon_2}$,
the rescaling formula (Lemma~\ref{lem:rescaling}) allows to rescale 
$k$ and $\ell$
by $\lambda \in \F_q^\x$, and  Proposition~\ref{prop:relations}
shows how to compute all
the values $\prod \calS_{k,\ell}^{++}$, $\prod \calS_{k,\ell}^{+-}$, 
$\prod \calS_{k,\ell}^{-+}$, $\prod \calS_{k,\ell}^{--}$ if just one of them
is known. Thus, we have reduced to computing one example of
$\prod \calS_{k,\ell}^{\pm,\pm}$ for each ratio 
$\tau=-k/\ell\in\F_q\cup\{\infty\}$, $\tau\ne -1$.  
It turns out that we can achieve this,
and the technique uses Dickson polynomials.
Once the results have been obtained for $\calS_{k,\ell}^{\pm, \pm}$, 
they can be transferred
to results about $\calT_{-k,\ell}^{\pm,\pm}$ using (\ref{ST}).

Our policy is to normalize $k$ and $\ell$ for $\calS_{k,\ell}^{\pm,\pm}$ 
so that $\ell-k=4$. 
That is, given $\tau \ne -1$, we set $k=-4\tau/(\tau+1)$ and $\ell=4/(\tau+1)$.
We set $r=\ell-2$, so $k=r-2$, $\ell=r+2$, and
$\tau=(2-r)/(2+r)$. For purpose of applying Proposition~\ref{prop:relations},
it is useful to note that for $\tau\ne \infty$,
\begin{equation} \label{jacobikell} 
\jacobi k q = \varepsilon \jacobi \tau q \jacobi{\tau+1}q,\qquad
\jacobi \ell q = \jacobi{\tau+1}q,\qquad \jacobi{\ell-k}q=1.
\end{equation}
It is also convenient to set $\tau'=k/4$. Then $\tau'+1=\ell/4$.
When discussing $\calT_{j,\ell}^{\pm,\pm}$, we set $j=-k$, so $\ell+j=4$.

Because these normalizations arise frequently, we highlight the relations between
the variables $\tau,k,\ell,r,j,\tau'$. 
\begin{eqnarray}
\tau &\in& \F_q \cup \{\infty\},\qquad \tau \ne -1 \label{eq_tau}\\
k&=&-j=-4\tau/(\tau+1)=4\tau'=r-2 \in \F_q \label{eq_k} \\
\ell&=&4+k=4-j=4/(\tau+1)=4(\tau'+1) = r+2 \in \F_q  \label{eq_ell} \\
r &=& k+2=2-j=\ell-2 \label{eq_r} \\
\tau'&=& k/4=-\tau/(\tau+1),\qquad \tau'+1=\ell/4=1/(\tau+1) \label{eq_tauPrime} \\
\tau &=& j/\ell=-k/\ell = (2-r)/(2+r) = -\tau'/(\tau'+1) \label{eq_tau2} 
\end{eqnarray}

With these normalizations, our formulas for 
$\prod \calS_{k,\ell}^{\pm,\pm}$ and $\prod \calT_{j,\ell}^{\pm,\pm}$
are summarized in Theorems~\ref{thm:main1} and~\ref{thm:main2} and Tables~\ref{tableMain1}--\ref{tableTjl2}
below.  Theorem~\ref{thm:main1} and Tables~\ref{tableMain1},~\ref{tableTjl1} handle all cases where 
$\tau=-k/\ell \in\{0,\infty, 1,1/3,3\}$ or where $\tau$ and $\tau+1$ are squares.  
Theorem~\ref{thm:main2} and Tables~\ref{tableMain2},~\ref{tableTjl2} handle all cases
where $\tau$ and/or $\tau+1$ are nonsquares. In the tables,
the conditions on $\tau$ are given on the left, and the
corresponding values for normalized $\prod \calS_{k,\ell}^{\pm,\pm}$ or $\prod \calT_{j,\ell}^{\pm,\pm}$
are given in the columns on the right. In some instances, 
separate formulas must be given depending on $q \bmod 4$, $q\bmod 8$, or $q \bmod{12}$.

\begin{theorem} \label{thm:main1}
Let $\tau$, $j$, $k$, $\ell$ as in (\ref{eq_tau})--(\ref{eq_tau2}) and $\varepsilon = \jacobi{-1}q$. 
The following formulas hold for $\prod \calS_{k,\ell}^{\pm,\pm}$ and $\prod \calT_{j,\ell}^{\pm,\pm}$ 
when $\tau =j/\ell\in \{0,\infty,1,1/3,3\}\cup \calA_{0,1}^{++}$. 
(For each $\tau$, we provide a formula for one choice of $(\pm,\pm)$, as the remaining ones
can be derived from that one using methods of Section~\ref{sec:relations}.
For the reader's convenience, we provide in
Tables~\ref{tableMain1} and~\ref{tableTjl1} the
values of $\prod \calS_{k,\ell}^{\varepsilon_1,\varepsilon_2}$ and 
$\prod \calT_{j,\ell}^{\varepsilon_1,\varepsilon_2}$ 
for all $\varepsilon_1,\varepsilon_2$.)
\begin{enumerate}
\item[({\it i})] $\tau=0$ or $\infty$: \quad
$\prod \calS_{0,4}^{-\varepsilon,-} = \prod\calT_{0,4}^{--} =\jacobi 2 q \cdot 2$;\quad
$\prod \calS_{-4,0}^{-\varepsilon,-}=\prod\calT_{4,0}^{--} =2$.
\item[({\it ii})] $\tau = 1$: \quad  
If $q \equiv \pm1 \pmod 8$ then 
$\prod \calS_{-2,2}^{-\varepsilon,-} = \prod\calT_{2,2}^{--} 
= (-1)^{(q-\varepsilon)/8}\cdot 2$,
and
if $q \equiv \pm3 \pmod 8$ then $\prod \calS_{-2,2}^{\varepsilon,+} = \prod\calT_{2,2}^{++}   
= (-1)^{(q-\varepsilon-4)/8}=(-1)^{\lfloor q/8 \rfloor}$.
\item[({\it iii})] $\tau = 1/3$ or $3$:  \quad
$\prod \calS_{-1,3}^{-\varepsilon,-} = \prod\calT_{1,3}^{--}
= \begin{cases}2 & \text{if $q \equiv \pm1 \pmod{12}$} \\
-1& \text{otherwise} 
\end{cases} $  \\
and 
$\prod \calS_{-3,1}^{-\varepsilon,-} = \prod\calT_{3,1}^{--} = \jacobi 2q \cdot \prod \calT_{1,3}^{--}$.
\item[({\it iv})]
Suppose that $\tau \in \calA_{0,1}^{++}$.
(Equivalently, $j$ and $\ell$ are nonzero squares.) Then
$$\prod \calS_{-j,\ell}^{-\varepsilon,-} = \prod \calT_{j,\ell}^{--}
= \jacobi{4\pm 2 \sqrt{\ell}}q \cdot 2 = \jacobi{2\pm \sqrt j}q \cdot 2.$$
(The Legendre symbols in the above formula are well defined and equal to one another by~(\ref{jell}).)
Also, $\prod \calT_{\ell,j}^{--}=\jacobi 2 q \prod \calT_{j,\ell}^{--}$.
\end{enumerate}
\end{theorem}

Theorem~\ref{thm:main1} addresses all cases where $\tau$ and $\tau+1$
are squares.
If $\tau$ and/or $\tau+1$ is a nonsquare, then exactly one of
$\tau$, $\tau+1$, and $\tau/(\tau+1)$ is a square.
Accordingly, the next theorem gives a formula for
$\prod \calS_{k,\ell}^{\pm,\pm}$ in terms of a square root of
$\tau$, $\tau+1$, or $\tau/(\tau+1)$. This appears to give the formula
only up to a choice of sign, but in Section~\ref{subsec:canonical} 
a prescription is given for which square root to take.
For each $\tau$, we give $\prod \calS_{k,\ell}^{\varepsilon_1,\varepsilon_2}$ for
normalized $k,\ell$ and for one value of $\varepsilon_1$ and $\varepsilon_2$.
The others can be computed by methods of Section~\ref{sec:relations}.
For the reader's convenience, we provide
in Tables~\ref{tableMain2} and~\ref{tableTjl2}
the values of $\prod \calS_{k,\ell}^{\varepsilon_1,\varepsilon_2}$ and 
$\prod \calT_{-k,\ell}^{\varepsilon_1,\varepsilon_2}$ for normalized
$k$ and $\ell$ and 
for all $\varepsilon_1,\varepsilon_2$.

\begin{theorem} \label{thm:main2}  Let $\tau \in \F_q$, $\tau \ne -1$,
$j=-k=4\tau/(\tau+1)$, $\ell=4/(\tau+1)$.
\begin{enumerate}
\item[({\it i})] Suppose $\tau$ is a nonzero square and 
$\tau+1$ is a nonsquare.  (Equivalently, $j$ and $\ell$ are nonsquares.)
Then
$$\prod\calT_{j,\ell}^{++} = \prod \calS_{k,\ell}^{\varepsilon,+}
= -\jacobi 2 q\cdot \frac{\tau+1}{2\,\,\sqrt{\tau\,\,}}
= -\jacobi 2q \cdot \frac{2}{\ell\,\,\sqrt{j/\ell\,\,}}.$$
\item[({\it ii})] Suppose $\tau$ is a nonsquare and $\tau+1$ is a square.  
(Equivalently, $j$ is a nonsquare and $\ell$ is a square.)
Then
$$\prod\calT_{j,\ell}^{--} = \prod \calS_{k,\ell}^{-\varepsilon,-}
= \jacobi 2 q \cdot 2/\sqrt{\tau+1} = \jacobi 2q \sqrt{\ell\,}. $$
\item[({\it iii})] Suppose $\tau$ and $\tau+1$ are nonsquares.   
(Equivalently, $j$ is a square and $\ell$ is a nonsquare.)
Then
$$\prod\calT_{j,\ell}^{--} = \prod \calS_{k,\ell}^{-\varepsilon,-}
= 2 \sqrt{\tau/(\tau+1)}=\sqrt{j\,}.$$
\end{enumerate}
In these formulas, the square roots 
are determined by the prescription given in 
Section~\ref{subsec:canonical}.
\end{theorem}

Formula~(\ref{T22}) of the introduction can be found in the row ``$\tau=1$'' of Table~\ref{tableTjl1}.
Formula (\ref{T13}) follows from the row ``$\tau=1/3$'' of Table~\ref{tableTjl1}
if $q$ is not a power of~3, and from the row 
``$\tau=\infty$'' of Table~\ref{tableTjl1} if $q$ is a power of~3.
Formulas (\ref{sqrtell}) and (\ref{sqrtj}) follow from 
Theorem~\ref{thm:main2}({\it ii}) and ({\it iii}).
Formula~(\ref{sqrtellj}) follows from Theorem~\ref{thm:main2}({\it i})
and Proposition~\ref{prop:relations} (or see Table~\ref{tableTjl2}).

\begin{center} 
\begin{table}
\scalebox{0.85}[0.90]{
\begin{tabular}{|c|c|c|c|c|c|c|c|}
\hline
\raisebox{-0.5ex}{Condition} 
& \raisebox{-1.2ex}{$k$} & \raisebox{-1.2ex}{$\ell$} 
&\raisebox{-0.3ex}{Condition}
& \raisebox{-1.2ex}{$\prod \calS_{k,\ell}^{++}$} &
\raisebox{-1.2ex}{$\prod \calS_{k,\ell}^{+-}$} 
&\raisebox{-1.2ex}{$\prod \calS_{k,\ell}^{-+}$} 
& \raisebox{-1.2ex}{$\prod \calS_{k,\ell}^{--}$}  \\
\raisebox{0.5ex}{on $\tau$} &&& \raisebox{0.5ex}{on $q$} &&&& \\
\hline\hline 
\raisebox{-2.6ex}{$\tau=0$} &  \raisebox{-2.6ex}{0} & \raisebox{-2.6ex}{4} &$1 \pmod 4$ &  \rule[-4mm]{0mm}{10mm} $\jacobi2q/4$ & $\jacobi2q$ 
                               & $\jacobi2q/2$ & $\jacobi2q\cdot 2$ \\
\cline{4-8}
         &   &   & $3 \pmod 4$ &  \rule[-4mm]{0mm}{10mm} $\jacobi2q/2$ & $\jacobi2q\cdot 2$ 
                 & $-\jacobi2q/4$ & $-\jacobi2q$ \\
\hline  
\raisebox{-1.2ex}{$\tau=\infty$} &\raisebox{-1.2ex}{$-4$} & 
\raisebox{-1.2ex}{0} &  \rule[-1mm]{0mm}{5mm} $1 \pmod 4$ & $-1/4$ & $1/2$ & 1 & 2 \\
\cline{4-8}
&      &   & \rule[-1mm]{0mm}{5mm}  $3 \pmod 4$ & 1 & 2 & $1/4$ & $-1/2$ \\
\hline 
&      &   &  $1 \pmod 8$ & \rule[-2mm]{0mm}{7mm}  $(-1)^{(q-1)/8}/8$   & $(-1)^{(q-1)/8}$  
                     & $(-1)^{(q-1)/8}$  & $(-1)^{(q-1)/8}\cdot 2$ \\
\cline{4-8}
\raisebox{-1.2ex}{$\tau=1$} & \raisebox{-1.2ex}{$-2$} &\raisebox{-1.2ex}{2} &  $5 \pmod 8$ & \rule[-2mm]{0mm}{7mm}  $(-1)^{(q-5)/8}$ & $(-1)^{(q-5)/8}$  
                     & $(-1)^{(q+3)/8}$   & $(-1)^{(q+3)/8}/4$  \\
\cline{4-8}
&      &    &  $3 \pmod 8$ & \rule[-2mm]{0mm}{7mm}  $(-1)^{(q-3)/8}$  & $(-1)^{(q-3)/8}/4$  
                      & $(-1)^{(q-3)/8}$ & $(-1)^{(q-3)/8}$  \\
\cline{4-8}
&      &    &  $7 \pmod 8$ & \rule[-2mm]{0mm}{7mm}  $(-1)^{(q+1)/8}$  & $(-1)^{(q+1)/8}\cdot 2$ 
                      & $(-1)^{(q-7)/8}/8$  & $(-1)^{(q-7)/8}$ \\
\hline
&      &   & $1 \pmod {12}$ & \rule[-4mm]{0mm}{10mm} $\jacobi2q/6$ & $\jacobi2q$ 
                     & $\jacobi2q$ & $\jacobi2q\cdot 2$ \\
\cline{4-8}
\raisebox{-2.0ex}{$\tau=3$} & \raisebox{-2.0ex}{$-3$} & \raisebox{-2.0ex}{1}  
& $5 \pmod {12}$ & \rule[-4mm]{0mm}{10mm} $-\jacobi2q$  & $-\jacobi2q\cdot 2$  
                 & $-\jacobi2q/6$  & $-\jacobi2q$ \\
\cline{4-8}
&     &   &$7 \pmod {12}$ & \rule[-4mm]{0mm}{10mm} $-\jacobi2q/6$ &  $-\jacobi2q$  
                    & $\jacobi2q$ & $\jacobi2q\cdot 2$  \\
\cline{4-8}
&     &   &$11 \pmod {12}$ & \rule[-4mm]{0mm}{10mm} $\jacobi2q$  & $\jacobi2q\cdot 2$ 
                    & $-\jacobi2q/6$    & $-\jacobi2q$   \\
\hline
&      &   &  $1 \pmod {12}$ & \rule[-1mm]{0mm}{5mm}  $1/6$ & 1 & 1 & 2 \\
\cline{4-8}
\raisebox{-1.0ex}{$\tau=1/3$} & \raisebox{-1.0ex}{$-1$} & \raisebox{-1.0ex}{3} 
&  $5 \pmod {12}$ & \rule[-1mm]{0mm}{6mm}  1 & $1/6$  & $-2$ & $-1$ \\
\cline{4-8}
&      &   &  $7 \pmod {12}$ & \rule[-1mm]{0mm}{5mm}  $-2$ & $-1$ & $-1$ & $-1/6$ \\
\cline{4-8}
&      &   &  $11 \pmod {12}$ & \rule[-1mm]{0mm}{5mm}  1 & 2 & $-1/6$ & $-1$ \\
\hline
\rule[-6mm]{0mm}{8mm} \raisebox{-1.0ex}[0pt]{$\tau$, $\tau+1$, and} & & & \raisebox{-1.5ex}[0pt]{$1 \pmod 4$} 
    & \raisebox{-1.5ex}[0pt]{$-1/(2k\ell)$}
    & \raisebox{-1.5ex}[0pt]{1} 
    & \raisebox{-1.5ex}[0pt]{1}  & 
      \raisebox{-1.5ex}[0pt]{2} \\
\cline{4-8}
\raisebox{0.5ex}[0pt]{$1\pm 1/\sqrt{\tau+1}$} & \raisebox{0.8ex}{$\frac{-4\tau}{\tau+1}$}
& \raisebox{0.8ex}{$\frac{4}{\tau+1}$} & &&&& \\
\raisebox{0.5ex}[0pt]{are squares}&  &
  &\raisebox{1.5ex}[0pt]{ $ 3\pmod 4$ }
  &\raisebox{1.5ex}[0pt]{ 1}  
  &\raisebox{1.5ex}[0pt]{ 2 }
  &\raisebox{1.5ex}[0pt]{ $1/(2k\ell)$ }
   &\raisebox{1.5ex}[0pt]{ $-1$} \\
\hline
$\tau$, $\tau+1$ are & & 
& \raisebox{-1.5ex}[0pt]{$1 \pmod 4$ }
& \raisebox{-1.5ex}[0pt]{$1/(2k\ell)$}
& \raisebox{-1.5ex}[0pt]{$-1$} 
& \raisebox{-1.5ex}[0pt]{$-1$}  
& \raisebox{-1.5ex}[0pt]{$-2$} \\
\raisebox{0.5ex}[0pt]{squares;}  & & &             &  &  &  &  \\
\cline{4-8}
\raisebox{-0.5ex}[0pt]{$1\pm1/\sqrt{\tau+1}$ }
  &\raisebox{1.0ex}[0pt]{$\frac{-4\tau}{\tau+1}$}
  &\raisebox{1.0ex}[0pt]{ $\frac{4}{\tau+1}$} 
  & \raisebox{-1.5ex}[0pt]{$3 \pmod 4$ }
  &\raisebox{-1.5ex}[0pt]{ $-1$}  
  &\raisebox{-1.5ex}[0pt]{ $-2$ }
  &\raisebox{-1.5ex}[0pt]{ $(-1)/(2k\ell)$ }
  &\raisebox{-1.5ex}[0pt]{ 1} \\
are nonsquares&  & & &   &   &   &   \\
\hline
\end{tabular}
}
\caption{ $\prod \calS_{k,\ell}^{\pm,\pm}$ for values of 
$\tau=-k/\ell$ given
in Theorem~\ref{thm:main1}. Here $k$ and $\ell$ are normalized so
that $\ell-k=4$. One entry from each row is proved in Theorem~\ref{thm:main1},
and the remaining entries of the row can be deduced using methods of 
Section~\ref{sec:relations}.} In the bottom two rows, we note that
$1\pm 1/\sqrt{\tau+1}$ are squares iff $2\pm2\sqrt{\tau/(\tau+1)}$ 
are squares by Proposition~\ref{prop:consequence2}. 
Also, $1\pm1/\sqrt{\tau+1}=1\pm\sqrt{\ell\,\,}/2$
and $2\pm2\sqrt{\tau/(\tau+1)}=2\pm\sqrt{-k\,}$.
\label{tableMain1}
\end{table}
\end{center}

\begin{center}
\begin{table}
\scalebox{0.85}[0.90]{
\begin{tabular}{|c|c|c|c|c|c|c|}
\hline  
Condition 
& \raisebox{-1.0ex}{\parbox[c][0pt]{8pt}{$j$}}
& \raisebox{-1.0ex}{\parbox[c][0pt]{8pt}{$\ell$}}
& \raisebox{-1.0ex}{\parbox[c][0pt]{60pt}{\quad$\prod \calT_{j,\ell}^{++}$} }
& \raisebox{-1.0ex}{\parbox[c][0pt]{60pt}{\quad$\prod \calT_{j,\ell}^{+-}$}}
& \raisebox{-1.0ex}{\parbox[c][0pt]{60pt}{\quad$\prod \calT_{j,\ell}^{-+}$}}
& \raisebox{-1.0ex}{\parbox[c][0pt]{60pt}{\quad$\prod \calT_{j,\ell}^{--}$}} \\
on $\tau$ and $q$ &&& &&& \\
\hline\hline  
$\tau=0$ &  0 & 4 &    \rule[-4mm]{0mm}{10mm} $\jacobi{-2}q/4$ & $\jacobi{-2}q$  &
                       $\jacobi2q/2$ & $\jacobi2q\cdot 2$ \\
\hline
$\tau=\infty$ & 4 & 0  & \rule[-4mm]{0mm}{10mm} $- \jacobi{-1}q /4$ & $\jacobi{-1}q/2$ & 1 & 2 \\
\hline  
\raisebox{-0.3ex}[0pt]{$\tau = 1$}
& \raisebox{-1.5ex}[0pt]{\parbox[c][0pt]{8pt}{\,\,\,\, 2}} 
& \raisebox{-1.5ex}[0pt]{\parbox[c][0pt]{8pt}{\,\,\,\, 2}}
& \raisebox{-2.0ex}[0pt]{\parbox[c][0pt]{60pt}{$(-1)^{\lfloor q/8\rfloor}/8$}}
& \raisebox{-2.0ex}[0pt]{\parbox[c][0pt]{60pt}{\quad$(-1)^{\lfloor q/8 \rfloor}$  }}
& \raisebox{-2.0ex}[0pt]{\parbox[c][0pt]{60pt}{
\raisebox{0.3ex}{$(-1)^{\lfloor(q+3)/8\rfloor}$}  }}
& \raisebox{-2.0ex}[0pt]{\parbox[c][0pt]{80pt}{
 $(-1)^{\lfloor (q+3)/8\rfloor}\cdot 2$}}  \\
\raisebox{0.3ex}{$q \equiv \pm1 \pmod 8$} & & &&&& \\
\hline
\raisebox{-0.3ex}[0pt]{$\tau = 1$}
& \raisebox{-1.5ex}[0pt]{\parbox[c][0pt]{8pt}{\,\,\,\, 2}} 
& \raisebox{-1.5ex}[0pt]{\parbox[c][0pt]{8pt}{\,\,\,\, 2}}
& \raisebox{-2.0ex}[0pt]{\parbox[c][0pt]{40pt}{
$(-1)^{\lfloor q/8\rfloor}$ }} 
& \raisebox{-2.0ex}[0pt]{\parbox[c][0pt]{40pt}{   
$(-1)^{\lfloor q/8\rfloor}$  }}
& \raisebox{-2.0ex}[0pt]{\parbox[c][0pt]{60pt}{
$(-1)^{\lfloor(q+3)/8\rfloor}$   }}
&\raisebox{-2.0ex}[0pt]{\parbox[c][0pt]{80pt}{$(-1)^{\lfloor(q+3)/8\rfloor}/4$}}  \\
$q = \pm 3 \pmod 8$ & & &&&& \\
\hline  
\raisebox{-0.5ex}[0pt]{\parbox[c][0pt]{30pt}{$\tau = 3$ }}
& \raisebox{-1.5ex}[0pt]{\parbox[c][0pt]{8pt}{\,\,\,\, 3}} 
& \raisebox{-1.5ex}[0pt]{\parbox[c][0pt]{8pt}{\,\,\,\, 1}}
& \raisebox{-2.0ex}{\parbox[c][0pt]{60pt}{\quad$\jacobi{-2}q/6$}}
& \raisebox{-2.0ex}{\parbox[c][0pt]{60pt}{\ \quad$\jacobi{-2}q$}}
& \raisebox{-2.0ex}{\parbox[c][0pt]{60pt}{\ \quad$\jacobi{2}q$}} 
& \raisebox{-2.0ex}{\parbox[c][0pt]{80pt}{\ \quad$\jacobi{2}q\cdot 2$}} \\
\raisebox{0.3ex}{$q \equiv \pm 1 \pmod{12}$} &&&&&& \\
\hline
\raisebox{-0.5ex}[0pt]{\parbox[c][0pt]{30pt}{$\tau = 3$ }}
& \raisebox{-1.5ex}[0pt]{\parbox[c][0pt]{8pt}{\,\,\,\, 3}} 
& \raisebox{-1.5ex}[0pt]{\parbox[c][0pt]{8pt}{\,\,\,\, 1}}
& \raisebox{-2.0ex}{\parbox[c][0pt]{60pt}{\quad$-\jacobi{-2}q$}}
& \raisebox{-2.0ex}{\parbox[c][0pt]{60pt}{\ $-\jacobi{-2}q \cdot 2$  }}
& \raisebox{-2.0ex}{\parbox[c][0pt]{60pt}{\quad $-\jacobi{2}q /6$  }}
& \raisebox{-2.0ex}{\parbox[c][0pt]{60pt}{\quad $-\jacobi{2}q$}} \\
$q \equiv \pm5 \pmod {12}$ & &&&&& \\
\hline
\raisebox{-0.5ex}[0pt]{\parbox[c][0pt]{50pt}{$\tau = 1/3$ }}
& \raisebox{-1.5ex}[0pt]{\parbox[c][0pt]{8pt}{\,\,\,\, 1}}
& \raisebox{-1.5ex}[0pt]{\parbox[c][0pt]{8pt}{\,\,\,\, 3}} 
& \raisebox{-2.0ex}{\parbox[c][0pt]{60pt}{\quad$\jacobi{-1}q/6$ }}
& \raisebox{-2.0ex}{\parbox[c][0pt]{60pt}{\quad$\jacobi{-1}q$}} 
& \raisebox{-2.0ex}{\parbox[c][0pt]{60pt}{\qquad1}}
& \raisebox{-2.0ex}{\parbox[c][0pt]{60pt}{\qquad2}} \\
$q\equiv \pm 1 \pmod {12}$ &&&&&& \\
\hline  
\raisebox{-0.5ex}[0pt]{\parbox[c][0pt]{50pt}{$\tau = 1/3$ }}
& \raisebox{-1.5ex}[0pt]{\parbox[c][0pt]{8pt}{\,\,\,\, 1}}
& \raisebox{-1.5ex}[0pt]{\parbox[c][0pt]{8pt}{\,\,\,\, 3}} 
& \raisebox{-2.0ex}{\parbox[c][0pt]{60pt}{\quad$\jacobi{-1}q$ }}
& \raisebox{-2.0ex}{\parbox[c][0pt]{60pt}{\quad$\jacobi{-1}q/6$ }}
& \raisebox{-2.0ex}{\parbox[c][0pt]{60pt}{\qquad$-2$}}
& \raisebox{-2.0ex}{\parbox[c][0pt]{60pt}{\qquad$-1$}} \\
$q \equiv \pm5 \pmod {12}$ &&& &&& \\
\hline  
\rule[-6mm]{0mm}{8mm} \raisebox{-1.0ex}[0pt]{$\tau$, $\tau+1$, and} & & & & & & \\
\raisebox{0.5ex}[0pt]{$1\pm1/\sqrt{\tau+1}$} & $\frac{4\tau}{\tau+1}$  & $\frac{4}{\tau+1}$ 
    & \raisebox{0.5ex}[0pt]{\rule[-4mm]{0mm}{10mm} $\jacobi{-1}q/(2j\ell)$} 
    & \raisebox{0.5ex}[0pt]{$\jacobi{-1}q$} 
    & \raisebox{0.5ex}[0pt]{$1$}  & 
      \raisebox{0.5ex}[0pt]{2} \\
\raisebox{0.5ex}[0pt]{are nonzero squares}&  & & & & & \\
\hline
$\tau$, $\tau+1$ are & & &           &  &  &   \\
\raisebox{0.5ex}[0pt]{nonzero squares;}   && & & & & \\
\raisebox{-0.5ex}[0pt]{$1\pm 1/\sqrt{\tau+1}$ } 
  &\raisebox{0.5ex}[0pt]{$\frac{4\tau}{\tau+1}$}
  &\raisebox{0.5ex}[0pt]{ $\frac{4}{\tau+1}$}  
    & \raisebox{0.5ex}[0pt]{\rule[-4mm]{0mm}{10mm} $-\jacobi{-1}q/(2j\ell)$} 
    & \raisebox{0.5ex}[0pt]{$-\jacobi{-1}q$} 
    & \raisebox{0.5ex}[0pt]{$-1$}  & 
      \raisebox{0.5ex}[0pt]{$-2$} \\
are nonsquares&  & & &   &   &      \\
\hline
\end{tabular}
}
\caption{ $\prod \calT_{j,\ell}^{\pm,\pm}$ for values of 
$\tau=j/\ell$ given in Theorem~\ref{thm:main1}.
Here $j$ and $\ell$ are normalized so that $j+\ell=4$, thus $j = 4\tau/(\tau+1)$ and
$\ell=4/(\tau+1)$. In the bottom two rows of the table, note that
$1 \pm 1/\sqrt{\tau+1}$ is a square iff $2\pm 2\sqrt{\tau/(\tau+1)}$
is a square by Proposition~\ref{prop:consequence2}. Also note that
$1\pm 1/\sqrt{\tau+1}=1 \pm \sqrt{\ell\,}/2$ and 
$2\pm 2\sqrt{\tau/(\tau+1)}=2\pm \sqrt j$.}
\label{tableTjl1}
\end{table}
\end{center}

\begin{center}
\begin{table}
\begin{tabular}{|c|c|c|c|c|c|c|}
\hline 
\raisebox{-0.7ex}{Condition} &\raisebox{-1.5ex}{Value of $c$} 
&\raisebox{-0.7ex}{Condition} & \raisebox{-1.5ex}{$\prod \calS_{k,\ell}^{++}$} &
\raisebox{-1.5ex}{$\prod \calS_{k,\ell}^{+-}$}
& \raisebox{-1.5ex}{$\prod \calS_{k,\ell}^{-+}$} 
& \raisebox{-1.5ex}{$\prod \calS_{k,\ell}^{--}$}  \\
\raisebox{0.7ex}{on $\tau$} & &\raisebox{0.7ex}{on $q$} &&&& \\
\hline\hline
\rule[-4mm]{0mm}{10mm} $\jacobi{\tau}q=1$, & 
\raisebox{-2.0ex}[0pt]{\parbox[c][0pt]{40pt}{$\jacobi2q \sqrt{\tau}$}} 
& $1\pmod 4$ 
& $\frac{-(\tau+1)}{2c}$ & $-c$ & $1/c$ & $\frac{\tau+1}{8c}$ \\
\cline{3-7}
\rule[-4mm]{0mm}{10mm} $\jacobi{\tau+1}q=-1$ & & $3\pmod 4$ & $-1/c$ & $\frac{-(\tau+1)}{8c}$
   &$\frac{-(\tau+1)}{2c}$ & $-c$ \\
\hline
\rule[-4mm]{0mm}{10mm} $\jacobi{\tau}q=-1$, 
&\raisebox{-2.0ex}[0pt]{\parbox[c][0pt]{60pt}{$\jacobi2q\sqrt{\tau+1}$ }}
& $1\pmod 4$ & $c/2$ & $c$
   & $\frac{c(\tau+1)}{16\tau}$ & $2/c$  \\
\cline{3-7}
\rule[-4mm]{0mm}{10mm} $\jacobi{\tau+1}q=1$ & & $3\pmod 4$ &$\frac{c(\tau+1)}{16\tau}$ &$2/c$ 
  & $-c/2$&  $-c$\\
\hline
\rule[-4mm]{0mm}{10mm} $\jacobi{\tau}q=-1$, 
& \raisebox{-2.0ex}[0pt]{\parbox[c][0pt]{60pt}{$\sqrt{\frac{\tau}{\tau+1}}$}}
& $1\pmod 4$ & $-1/(2c)$ & $\frac{-(\tau+1)}{16c}$ & $1/c$ & $2c$  \\
\cline{3-7}
\rule[-4mm]{0mm}{10mm} $\jacobi{\tau+1}q=-1$ & & $3\pmod 4$ 
& $1/c$ &$2c$ & $1/(2c)$ & $\frac{\tau+1}{16c}$ \\
\hline
\end{tabular}
\caption{$\prod \calS_{k,\ell}^{\pm,\pm}$ for values of $\tau=-k/\ell$ given
in Theorem~\ref{thm:main2}. Here 
$k=-4\tau/(\tau+1)$ and $\ell=4/(\tau+1)$.
In each case,
$\tau \in \F_q \setminus \{0,-1\}$, and exactly one of $\tau$, $\tau+1$,
or $\tau/(\tau+1)$ is a square. Accordingly, we let $c\in\F_q$ denote
a square root of $\tau$, $\tau+1$, or $\tau/(\tau+1)$.
The choice for $\sqrt\tau$, $\sqrt{\tau+1}$, or $\sqrt{\tau/(\tau+1)}$ is
prescribed in Section~\ref{subsec:canonical}.}
\label{tableMain2}
\end{table}
\end{center}

\begin{center}
\begin{table}
\begin{tabular}{|c|c|c|c|c|c|c|}
\hline
\rule[-3mm]{0mm}{9mm} Condition on $\tau$ &$c$ & $\prod \calT_{j,\ell}^{++}$ &
$\prod \calT_{j,\ell}^{+-}$ & $\prod \calT_{j,\ell}^{-+}$ & $\prod \calT_{j,\ell}^{--}$  \\
\hline\hline
\rule[-4mm]{0mm}{10mm} $\jacobi{\tau}q=1$, $\jacobi{\tau+1}q=-1$ 
& $\jacobi2q \sqrt{\tau}$ 
& $\frac{-(\tau+1)}{2c}$ & $-c$ & $\varepsilon/ c$ 
& \raisebox{-0.2ex}[0pt]{\parbox[c][0pt]{50pt}{\quad$\frac{\varepsilon(\tau+1)}{8c}$}} \\
\hline
\rule[-4mm]{0mm}{10mm} $\jacobi{\tau}q=-1$, $\jacobi{\tau+1}q=1$ 
& $\jacobi2q\sqrt{\tau+1}$ & $\varepsilon c/2$ & $\varepsilon c$
   & $\frac{c(\tau+1)}{16\tau}$ & $2/c$  \\
\hline
\rule[-4mm]{0mm}{10mm} $\jacobi{\tau}q=\jacobi{\tau+1}q=-1$ & \raisebox{0ex}[0pt]{\parbox[c][0pt]{50pt}{$\sqrt{\frac{\tau}{\tau+1}}$}}
& $-\varepsilon/(2c)$ & 
$\frac{\varepsilon(\tau+1)}{-16c}$
& $1/c$ & $2c$  \\
\hline
\end{tabular}
\caption{$\prod \calT_{j,\ell}^{\pm,\pm}$ for values of $\tau = j/\ell$ given in Theorem~\ref{thm:main2}.
Here $j=4\tau/(\tau+1)$, $\ell=4/(\tau+1)$,
and $\varepsilon=\jacobi{-1}q=(-1)^{(q-1)/2}$.
In each case, $\tau\in\F_q \setminus \{0,-1\}$ and exactly one of $\tau$, $\tau+1$,
or $\tau/(\tau+1)$ is a square. Accordingly, we let $c\in\F_q$ denote
a square root of $\tau$, $\tau+1$, or $\tau/(\tau+1)$.
The choice for $\sqrt\tau$, $\sqrt{\tau+1}$, or $\sqrt{\tau/(\tau+1)}$ is
prescribed in Section~\ref{subsec:canonical}.}
\label{tableTjl2}
\end{table}
\end{center}

\vfill\eject

\section{Proofs of Theorems~\ref{thm:main1} and~\ref{thm:main2}} \label{sec:proofs}

Our main tool to prove formulas for $\prod \calS_{k,\ell}^{\pm,\pm}$ is
a result from \cite{Bluher-dickson} about Dickson polynomials.
Section~\ref{subsec:Dickson} provides a brief background on Dickson 
polynomials and recalls relevant results from \cite{Bluher-dickson}.
Section~\ref{subsec:main1Proof} uses
those results to prove Theorem~\ref{thm:main1}({\it i})--({\it iii}).
Section~\ref{subsec:rho_ru}  relates $\tau$ with two other parameters 
$r$ and $u$.
Section~\ref{subsec:canonical} describes a way to uniquely determine 
square roots of $\tau$, $\tau+1$,
or $\tau/(\tau+1)$ in the case where exactly one of these three is a square.
Then in Section~\ref{subsec:main2Proof} we prove 
Theorem~\ref{thm:main1}({\it iv}) and
Theorem~\ref{thm:main2}.

\subsection{Using Dickson polynomials to compute $\prod \calS_{k,\ell}^{\pm,\pm}$}
\label{subsec:Dickson}
\parbox[c][20pt]{30pt}{ $\phantom{hello world}$ }

Dickson polynomials, introduced in Leonard Eugene Dickson's PhD thesis 
in 1896, are closely related to Chebyshev polynomials and have been widely
studied; see \cite{Lidl}.
Dickson polynomials of the first kind, $D_k(x)\in\Z[x]$, are determined from the recursion 
\begin{equation} \label{Ddef} \text{$D_0(x)=2$, $D_1(x)=x$, and
$D_{k+2}(x) = x D_{k+1}(x) - D_k(x)$.} \end{equation}
Dickson polynomials of the second kind, $E_k(x) \in \Z[x]$, are
determined from the recursion
\begin{equation} \label{Edef} \text{$E_0(x)=1$, $E_1(x)=x$, and
$E_{k+2}(x) = x E_{k+1}(x) - E_k(x)$.} \end{equation}
It is well known, and can easily be shown by induction, that 
$D_k$ and $E_k$ satisfy the functional
equations
\begin{equation} \label{DEfunctional} D_k(\ang u) = \ang{u^k},\qquad
E_{k-1}(\ang u) = \frac{u^k - 1/u^k}{u-1/u} \end{equation}
where we use the notation $$\ang u = u + 1/u.$$

Define functions $f^{\varepsilon_1,\varepsilon_2}(x)\in\F_q[x]$ by
\begin{equation} \label{fij} f^{\varepsilon_1,\varepsilon_2}(x) 
= \prod\left\{(x-b) : b \in \calA_{-2,2}^{\varepsilon_1,\varepsilon_2} \right\}. 
\end{equation}

\begin{theorem} (\cite{Bluher-dickson})
 \label{thm:fij}  Let $\varepsilon=\jacobi{-1}q$, so
$\varepsilon\in\{+1,-1\}$ and 
$q\equiv \varepsilon \pmod 4$. In $\F_q[x]$ we have
$$ D_{(q-\varepsilon)/4}(x) = f^{-\varepsilon,-}(x),\qquad
E_{(q-\varepsilon-4)/4}(x) = f^{\varepsilon,+}(x).$$
\end{theorem}
\begin{proof} Although this was proved in \cite{Bluher-dickson}, we give
here a slightly different proof that is based on 
Corollary~\ref{cor:v2}. Let $m=(q-\varepsilon)/4$.

{\it (i)}\   $|\calA_{-2,2}^{-\varepsilon,-}|
=m$ by (\ref{AijCard2}), so $D_{m}(x)$ and    
$f^{-\varepsilon,-}(x)$ are monic polynomials of the same degree.
Thus, it suffices to prove that $D_m$ vanishes at all $b \in \calA_{-2,2}^{-\varepsilon,-}$.
By Corollary~\ref{cor:v2}, $b=\ang{v^2}$ where $v^{q-\varepsilon}=-1$
and $v^4\ne1$. Then $i=v^{(q-\varepsilon)/2}$ is a square root of $-1$, and
$$D_m(b) = D_m(\ang{v^2})=\ang{v^{2m}}=\ang{v^{(q-\varepsilon)/2}}
=\ang{i}=0.$$

{\it (ii)}\  By (\ref{AijCard2}), $|\calA_{-2,2}^{\varepsilon,+}|=m-1$,
so $E_{m-1}(x)$ and $f^{\varepsilon,+}(x)$ are monic polynomials of the
same degree. To show $E_{m-1}=f^{\varepsilon,+}$, it suffices to show 
$E_{m-1}(b)=0$ for all 
$b\in \calA_{-2,2}^{\varepsilon,+}$. By Corollary~\ref{cor:v2}, $b=\ang{v^2}$,
where $v^{q-\varepsilon}=1$ and $v^4\ne1$.  Then 
$$E_{m-1}(b)=E_{m-1}(\ang{v^2})
=\left( (v^2)^m - (v^2)^{-m} \right)/(v^2-v^{-2}).$$
Here the denominator does not vanish, since $v^4\ne1$.
On the other hand $v^{2m}\in \{1,-1\}$ since its square is~1, so the 
numerator does vanish, and $E_{m-1}(b)=0$.
\end{proof}

\medskip

\noindent{\bf Example 1.} Consider $q=13$.
The nonzero squares in $\F_{13}$ are $1,3,4,9,10,12$ and the nonsquares
are $2,5,6,7,8,11$. Thus, $\calA_{-2,2}^{++} = \{1,12\}$, $\calA_{-2,2}^{+-}=\{3,5,6\}$, 
$\calA_{-2,2}^{-+}=\{7,8,10\}$,
$\calA_{-2,2}^{--} = \{4,9,0\}$. 
We have $f^{++}(x)=(x-1)(x-12)\equiv x^2-1=E_2(x) \pmod{13}$ and 
$f^{--}(x)=(x-4)(x-9)(x-0) \equiv x^3-3x =D_3(x)\pmod{13}$. 

\medskip

\noindent{\bf Example 2.} Consider $q=23$.
The nonzero squares in $\F_{23}$ are $1,2,3,4,6,8,9,12,13,16,18$ and the nonsquares
are $5,7,10,11,14,15,17,19,20,21,22$. Thus, 
$$\calA_{-2,2}^{++} = \{4,6,10,11,14\},\qquad 
\calA_{-2,2}^{+-}=\{3,5,8,15,18,20\},$$
$$\calA_{-2,2}^{-+}=\{0,1,7,16,22\},\qquad \calA_{-2,2}^{--} = \{9,12,13,17,19\}.$$ 
Then
$f^{-+}(x)=x(x-1)(x-7)(x-16)(x-22)\equiv x^5-4 x^3+3x=E_5(x) \pmod{23},$
$$f^{+-}(x)=(x-3)(x-5)(x-8)(x-15)(x-18)(x-20) \equiv x^6-6 x^4 + 9 x^2 - 2 =D_6(x)\pmod{23}.$$ 

\medskip
\begin{lemma} \label{lem:fijSkl} 
If $r \in \F_q \setminus \calA_{-2,2}^{\varepsilon_1,\varepsilon_2}$ then 
$(-1)^{|\calA_{-2,2}^{\varepsilon_1,\varepsilon_2}|} 
f^{\varepsilon_1,\varepsilon_2}(r) 
= \prod \calS_{r-2,r+2}^{\varepsilon_1,\varepsilon_2}$.
\end{lemma}
\begin{proof} The set of roots of $f^{\varepsilon_1,\varepsilon_2}$ is precisely 
$\calA_{-2,2}^{\varepsilon_1,\varepsilon_2}$,
so the hypothesis that $r \not \in \calA_{-2,2}^{\varepsilon_1,\varepsilon_2}$ implies that $f^{\varepsilon_1,\varepsilon_2}(r)\ne 0$.
We have 
$$f^{\varepsilon_1,\varepsilon_2}(r) = 
\prod \left\{r-b : b \in \calA_{-2,2}^{\varepsilon_1,\varepsilon_2} \right\},$$
\begin{eqnarray*}
(-1)^{|\calA_{-2,2}^{\varepsilon_1,\varepsilon_2}|} 
f^{\varepsilon_1,\varepsilon_2}(r) &=& 
\prod \left\{b-r : b \in \calA_{-2,2}^{\varepsilon_1,\varepsilon_2} \right\}\\
&=& \prod \left\{ a : a + r \in \calA_{-2,2}^{\varepsilon_1,\varepsilon_2} 
\right\} \\
&=& \prod \left\{ a \in \F_q : \jacobi {a+r-2} q= \varepsilon_1,\ 
\jacobi {a+r+2}q=\varepsilon_2 \right\}.
\end{eqnarray*}
Since this is nonzero, 0 does not occur in the product, and so this equals
$$ \prod \left\{ a \in \F_q^\x : \jacobi {a+r-2} q= \varepsilon_1,\ 
\jacobi {a+r+2}q=\varepsilon_2 \right\} = \prod \calS_{r-2,r+2}^{\varepsilon_1,\varepsilon_2}$$
as claimed. 
\end{proof}

\begin{proposition} \label{prop:Sklu} 
\rule[-3mm]{0mm}{9mm} 
Let $r \in \F_q$ and write $r = \ang u$, where $u \in \cj \F_q^\x$.  
Let $\varepsilon =\jacobi{-1}q$. \\
\noindent({\it i})\ \rule[-3mm]{0mm}{9mm} 
If $r \not \in \calA_{-2,2}^{-\varepsilon,-}$, 
then $\prod \calS_{r-2,r+2}^{-\varepsilon,-}=\ang{(-u)^{(q-\varepsilon)/4}}
= \jacobi 2 q \ang {u^{(q-\varepsilon)/4}}$. \\
\noindent({\it ii})\ \rule[-3mm]{0mm}{9mm} 
If $r \not \in \calA_{-2,2}^{\varepsilon,+} \cup \{2,-2\}$, then 
\begin{eqnarray*}
\prod \calS_{r-2,r+2}^{\varepsilon,+} &=&
-((-u)^{(q-\varepsilon)/4}-(-u)^{-(q-\varepsilon)/4})/(u-1/u) \\
&=&
-\jacobi 2q (u^{(q-\varepsilon)/4}-u^{-(q-\varepsilon)/4})/(u-1/u). 
\end{eqnarray*}
 \\
\end{proposition}
\begin{proof} ({\it i}) Suppose $r\in\F_q\setminus \calA_{-2,2}^{-\varepsilon,-}$. Then
\begin{eqnarray*}
\prod \calS_{r-2,r+2}^{-\varepsilon,-} &=& 
(-1)^{|\calA_{-2,2}^{-\varepsilon,-}|} f^{-\varepsilon,-}(r) \qquad
\text{by Lemma~\ref{lem:fijSkl}} \\
&=& (-1)^{(q-\varepsilon)/4} D_{(q-\varepsilon)/4}(r)  \qquad 
\text{by Theorem~\ref{thm:fij} and (\ref{AijCard2})} \\
&=& (-1)^{(q-\varepsilon)/4} \ang{ u^{(q-\varepsilon)/4 }} \qquad 
\text{by~(\ref{DEfunctional})}\\
&=& \ang{(-u)^{(q-\varepsilon)/4}} \qquad \text{since $-\ang v = \ang{-v} = -(v+1/v).$}
\end{eqnarray*}
By (\ref{jacobi2q}), we also have $\prod \calS_{r-2,r+2}^{-\varepsilon,-}=
\jacobi 2q \ang {u^{(q-\varepsilon)/4}}$. \\
({\it ii}) Suppose $r \not \in \calA_{-2,2}^{\varepsilon,+}\cup\{2,-2\}$.
Note that $r\not\in\{2,-2\}$ implies $u\not \in\{1,-1\}$, so in particular
$u-1/u\ne 0$. We have
\begin{eqnarray*}
\prod \calS_{r-2,r+2}^{\varepsilon,+} &=& 
(-1)^{|\calA_{-2,2}^{\varepsilon,+}|} f^{\varepsilon,+}(r)\qquad
\text{by Lemma~\ref{lem:fijSkl}} \\
&=& (-1)^{(q-\varepsilon-4)/4} E_{(q-\varepsilon-4)/4}(r)\qquad
\text{by Theorem~\ref{thm:fij} and (\ref{AijCard2})}\\
&=& (-1)^{(q-\varepsilon-4)/4} 
(u^{(q-\varepsilon)/4}-u^{-(q-\varepsilon)/4})/(u-1/u)\qquad
\text{by (\ref{DEfunctional})} \\
&=& -\left((-u)^{(q-\varepsilon)/4}-(-u)^{-(q-\varepsilon)/4}\right)/(u-1/u).
\end{eqnarray*}
The right side is also equal to $-\jacobi2q 
\left(u^{(q-\varepsilon)/4}-u^{-(q-\varepsilon)/4}\right)/(u-1/u)$ 
by (\ref{jacobi2q}).
\end{proof}

\subsection{Proof of Theorem~\ref{thm:main1}({\it i})--({\it iii}).}
\label{subsec:main1Proof}
\parbox[c][20pt]{30pt}{ $\phantom{hello world}$ }

We are ready to prove parts ({\it i})--({\it iii}) of
Theorem~\ref{thm:main1}. 

\noindent ({\it i})\ $\tau=0$:
Set $r=2=\ang{1}$. By Proposition~\ref{prop:Sklu}({\it i}), 
$\prod \calS_{0,4}^{-\varepsilon,-}=\jacobi 2q \ang 1 = \jacobi 2 q \cdot 2$.

\noindent $\tau=\infty$:
Set $r=-2=\ang{-1}$. By Proposition~\ref{prop:Sklu}({\it i}), 
$\prod \calS_{-4,0}^{-\varepsilon,-}=\ang 1 = 2$.

\noindent ({\it ii}) $\tau = 1$:  Set $r=0=\ang{i}$, where $i^2=-1$.

If $q\equiv \pm1 \pmod 8$ then $\jacobi 2 q = 1$ by Lemma~\ref{lem:simple}, 
so $0\not\in\calA_{-2,2}^{-\varepsilon,-}$. By Proposition~\ref{prop:Sklu}({\it i}), 
$$\prod \calS_{-2,2}^{-\varepsilon,-} = \ang{(-i)^{(q-\varepsilon)/4}}
=\ang{(-1)^{(q-\varepsilon)/8}}=(-1)^{(q-\varepsilon)/8} \cdot 2.$$

If $q\equiv \pm3 \pmod 8$ then $\jacobi 2 q = -1$, so $0 \not \in \calA_{-2,2}^{\varepsilon,+}$.
By Proposition~\ref{prop:Sklu}({\it ii}), 
$$\prod \calS_{-2,2}^{\varepsilon,+} = 
(i^{(q-\varepsilon)/4}-i^{-(q-\varepsilon)/4})/(i-1/i).$$
Write $q=8m+4+\varepsilon$. Since $i^{(q-\varepsilon)/4}=(-1)^m \cdot i$,  the right side is
$(-1)^m(i-1/i)/(i-1/i)=(-1)^m$.  Thus, $\prod \calS_{-2,2}^{\varepsilon,+} = (-1)^{(q-4-\varepsilon)/8}=
(-1)^{\lfloor q/8 \rfloor}$.

\noindent ({\it iii}) $\tau = 1/3$:  Set $r=1=\ang{-\omega}$, where 
$\omega$ is a primitive cube root of unity.
Then $\jacobi{r-2}q=\jacobi{-1}q=\varepsilon$, so 
$r \not \in\calA_{-2,2}^{-\varepsilon,-}$. By Proposition~\ref{prop:Sklu},
$\prod \calS_{-1,3}^{-\varepsilon,-}=\ang{\omega^m}$,  
where $m=(q-\varepsilon)/4$.  This equals 2 if $3|m$, or $-1$ if
$3\nmid m$. Here $3|m \iff q \equiv \pm1 \pmod{12}$.

\noindent $\tau = 3$:   $\prod \calS_{-3,1}^{-\varepsilon,-} 
= \prod \calT_{3,1}^{--} = \jacobi2q \prod \calT_{1,3}^{--}$
by (\ref{ST}) and Proposition~\ref{prop:Tellj}.

\noindent{\bf Remark.}\quad
If $q\equiv \pm 3 \pmod{12}$ ({\it i.e.} $q=3^n$), then $\calS_{-3,1}^{\varepsilon_1,\varepsilon_2}
=\calS_{0,4}^{\varepsilon_1,\varepsilon_2}$ and $\calT_{3,1}^{\varepsilon_1,\varepsilon_2}=
\calT_{0,4}^{\varepsilon_1,\varepsilon_2}$ because $0=3$ and $1=4$ in $\F_q$. For this reason,
we exclude the cases $q \equiv \pm3 \pmod{12}$ from Tables~\ref{tableMain1} and~\ref{tableTjl1}
when $\tau=1/3$ or $\tau=3$.

\subsection{Relation between $\tau$, $r$, and $u$.} \label{subsec:rho_ru}
\parbox[c][20pt]{30pt}{ $\phantom{hello world}$ }

Let $\tau \in \F_q \setminus \{0,-1\}$, $k,\ell,r,\tau'$ be as in equations~(\ref{eq_tau})--(\ref{eq_tau2}) and
$\varepsilon = \jacobi{-1}q$.
Note that $\tau'\in \F_q \setminus \{0,-1\}$.

\begin{lemma}  \label{lem:tauvw}
Let $v=\sqrt{\tau+1}+\sqrt\tau$ and 
$w = \sqrt{\tau'}+\sqrt{\tau'+1}$ (for any choice of square roots). Then $v^4\ne 1$, $w^4\ne 1$, 
and $r=\ang{w^2}$.
For $\varepsilon_1,\varepsilon_2 \in \{1,-1\}$ we have the following equivalence:
\begin{equation*} \tau \in \calA_{0,1}^{\varepsilon_1,\, \varepsilon_2}
\iff \tau' \in \calA_{0,1}^{\varepsilon\varepsilon_1\varepsilon_2,\, \varepsilon_2}
\iff r \in \calA_{-2,2}^{\varepsilon\varepsilon_1\varepsilon_2,\, \varepsilon_2} \iff
v^{q-\varepsilon_1\varepsilon_2}=\varepsilon_2 \iff w^{q-\varepsilon\varepsilon_1} = \varepsilon_2.
\end{equation*}
\end{lemma}
\begin{proof} As noted in (\ref{disjointUnion}), $\F_q\setminus\{0,-1\}$ is a disjoint union of the sets $\calA_{0,1}^{\varepsilon_1,\varepsilon_2}$.
By Theorem~\ref{thm:correspondence}, $\tau \in \calA_{0,1}^{\varepsilon_1,\varepsilon_2} $ iff $ 
v^{q-\varepsilon_1\varepsilon_2}=\varepsilon_2$ and $v^4\ne 1$.
By (\ref{eq_tauPrime}), $\tau \in \calA_{0,1}^{\varepsilon_1,\varepsilon_2} $ iff $ \tau' \in \calA_{0,1}^{\varepsilon
\varepsilon_1\varepsilon_2,\varepsilon_2}$.
By Theorem~\ref{thm:correspondence}, $\tau' \in \calA_{0,1}^{\varepsilon\varepsilon_1\varepsilon_2,\varepsilon_2} 
$ iff $ w^{q-\varepsilon\varepsilon_1}=\varepsilon_2$ and $w^4 \ne 1$. Since $r=4\tau'+2$, we have
$\jacobi{r-2}q = \jacobi{\tau'}q$ and $\jacobi{r+2}q=\jacobi{\tau'+1}q$, so
$\tau' \in \calA_{0,1}^{\varepsilon\varepsilon_1\varepsilon_2,\varepsilon_2} $ iff $ 
r \in \calA_{-2,2}^{\varepsilon\varepsilon_1\varepsilon_2,\varepsilon_2}$. 
It remains only to show that $r=\ang{w^2}$. Indeed,
$$\ang{w^2}= (\sqrt{\tau'+1}+\sqrt{\tau'})^2 + (\sqrt{\tau'+1}-\sqrt{\tau'})^2
=2(\tau'+1)+2\tau'=r.$$
\end{proof}

Let $u=w^2$, so that $r=\ang u$.
Table \ref{table:rho_ru} summarizes the relation between $\tau$, $r$ and $u$.

We can obtain an explicit relation between $v$ and $w$ as follows. Recall
$v=\sqrt{\tau+1}+\sqrt\tau$ and $w=\sqrt{\tau'+1}+\sqrt{\tau'}$. Changing
the choice of square roots in these formulas would change $v$ or $w$ to another
element in its orbit.  To fix matters, select $\sqrt{\tau+1}$,  $\sqrt\tau$, and
$\sqrt{\tau'}$ arbitrarily. By (\ref{eq_tauPrime}),  $\tau'+1=1/(\tau+1)$, so we may set
$\sqrt{\tau'+1}=1/\sqrt{\tau+1}$.  Let $i=\sqrt{\tau'}\sqrt{\tau+1}/\sqrt\tau$. 
By~(\ref{eq_tauPrime}), $i^2=\tau'(\tau+1)/\tau = -1$.
We have 
$$w = \sqrt{\tau'+1}+\sqrt{\tau'}=(1+i\sqrt\tau)/\sqrt{\tau+1}.$$
Using that $\sqrt\tau=(v-1/v)/2$ and $\sqrt{\tau+1}=(v+1/v)/2$, we obtain (for a particular choice of
$v,w$ within their orbits):
\begin{equation}
w = \frac{2+i(v-1/v)}{v+1/v}. \label{vw_relation}
\end{equation}

\begin{table}
\begin{center}
\begin{tabular}{|c|c|c|}
\hline
\rule[-2mm]{0mm}{8mm}
Condition on $\tau$ & Condition on $r$ & Condition on $u$ \\
\hline
$\tau = 0$ & $r=2$ & $u=1$ \\
\hline
$\tau = \infty$ & $r=-2$ & $u=-1$ \\
\hline
\rule[-3mm]{0mm}{8mm} 
$\tau \in \calA_{0,1}^{\varepsilon,+}$
& $ r\in \calA_{-2,2}^{++} $ &
$u^{(q-1)/2}=1$, $u^2\ne 1$ \\
\hline
\rule[-3mm]{0mm}{8mm} 
$\tau \in \calA_{0,1}^{\varepsilon,-}$
& $ r\in \calA_{-2,2}^{--} $ &
$u^{(q-1)/2}=-1$, $u^2\ne 1$ \\
\hline
\rule[-3mm]{0mm}{8mm} 
$\tau \in \calA_{0,1}^{-\varepsilon,+}$
& $ r\in \calA_{-2,2}^{-+} $ &
$u^{(q+1)/2}=1$, $u^2\ne 1$ \\
\hline
\rule[-3mm]{0mm}{8mm} 
$\tau \in \calA_{0,1}^{-\varepsilon,-}$
& $ r\in \calA_{-2,2}^{+-} $ &
$u^{(q+1)/2}=-1$, $u^2\ne 1$ \\
\hline
\rule[-3mm]{0mm}{8mm} 
$\tau \in \calA_{0,1}^{++}$
& $ r\in \calA_{-2,2}^{\varepsilon,+} $ &
$u^{(q-\varepsilon)/2}=1$, $u^2\ne 1$ \\
\hline
\rule[-3mm]{0mm}{8mm} 
$\tau \in \calA_{0,1}^{+-}$
& $ r\in \calA_{-2,2}^{-\varepsilon,-} $ &
$u^{(q-\varepsilon)/2}=-1$ \\
\hline
\rule[-3mm]{0mm}{8mm} 
$\tau \in \calA_{0,1}^{-+}$
& $ r\in \calA_{-2,2}^{-\varepsilon,+} $ &
$u^{(q+\varepsilon)/2}=1$, $u\ne 1$ \\
\hline
\rule[-3mm]{0mm}{8mm} 
$\tau \in \calA_{0,1}^{--}$
& $ r\in \calA_{-2,2}^{\varepsilon,-} $ &
$u^{(q+\varepsilon)/2}=-1$, $u\ne -1$ \\
\hline
\end{tabular}
\caption{Relation between $\tau$, $r$ and $u$ when $\tau \in\F_q\cup\infty$,
$\tau \ne -1$, $\tau = (2-r)/(2+r)$, $r=\ang u$. (When $\tau\in \{0,\infty\}$, use~(\ref{eq_k}), 
	and when $\tau\in\F_q \setminus \{0,-1\}$, use Corollary~\ref{cor:v2}.)}
\label{table:rho_ru}
\end{center}
\end{table}

If $\tau$ and $\tau+1$ are nonzero squares, more can be said.
In this case, Lemma~\ref{lem:tauvw} tells us that
$r \in \calA_{-2,2}^{\varepsilon,+}$ and $w^{q-\varepsilon}=1$.
In particular, $r+2$ and $-(r-2)$ are nonzero squares.  Let $\sqrt{r+2}$ and $\sqrt{2-r}$ denote
any choice of square roots.  
Note that $u^{(q-\varepsilon)/2}=w^{q-\varepsilon}=1$, and so $u^{(q-\varepsilon)/4} \in \{1,-1\}$.

\begin{proposition} \label{prop:ru2}
Suppose that $\tau \in \calA_{0,1}^{++}$, and let $r=\ang u\in\calA_{-2,2}^{\varep,+}$
be as in Table~\ref{table:rho_ru}. As per Proposition~\ref{prop:consequence2}, let
$c$ be a square root of $1+\tau$, $c'$ a square root of $1+1/\tau$, 
$$\nu = \jacobi{1\pm 1/c}q\qquad{\rm and}\qquad \mu = \jacobi{1\pm 1/c'}q.$$
Then
$$\nu = \jacobi 2q \mu,\qquad \sqrt{r+2} \in \calA_{-2,2}^{\varep\mu,\mu},
\qquad \sqrt{2-r} \in\calA_{-2,2}^{\varep\nu,\nu},\qquad u^{(q-\varep)/4}=\mu.$$
\end{proposition}
\begin{proof} 
By Proposition~\ref{prop:consequence2}, $\nu$ and $\mu$ are well defined
and $\nu = \jacobi 2q \mu$. Now $2/c\in\calA_{-2,2}^{\varep\mu,\mu}$ because
$$\jacobi {2/c+2}q=\jacobi{2(1+1/c)}q=\jacobi2q\nu=\mu,$$
$$\jacobi {2/c-2}q=\jacobi{-1}q \jacobi 2q\jacobi{1-1/c}q=\varepsilon\mu.$$
The same reasoning applied to $-c$ shows that 
$-2/c \in\calA_{-2,2}^{\varepsilon\mu,\mu}$ as well.
Since $(2/c)^2=4/(\tau+1)=r+2$, this shows
$\sqrt{r+2} = \pm 2/c \in \calA_{-2,2}^{\varepsilon\mu,\mu}$.

Now we show $u^{(q-\varep)/4}=\mu$. Write $r'=\sqrt{r+2}=\ang{s^2}$, 
where $s\in\cj\F_q$. Since 
$r'\in\calA_{-2,2}^{\varepsilon\mu,\mu}$, we have
$s^{q-\varepsilon}=\mu$ and $s^4\ne1$ by Corollary~\ref{cor:v2}.
Now $r+2=(r')^2=\ang{s^2}^2=\ang{s^4}+2$, so $r=\ang{s^4}$. Since also $r=\ang{u}$, 
this shows $u=s^4$ or $u=s^{-4}$. Thus, $s^{q-\varepsilon}=\mu$ implies 
that $u^{(q-\varep)/4}=s^{\pm(q-\varep)} = \mu^{\pm1}$.
Since $\mu=\mu^{-1}$, we conclude that $u^{(q-\varep)/4}=\mu$.

It remains only to prove that $\sqrt{2-r}\in \calA_{-2,2}^{\varep\nu,\nu}$.
Let $\widetilde\tau=1/\tau$, $\widetilde r=-r$, and $\widetilde u=-u$.  
Then $\widetilde\tau \in \calA_{0,1}^{++}$,
$\widetilde \tau=(2-\widetilde r)/(2+\widetilde r)$ and $\widetilde r=\ang{\widetilde u}$.  
Note that replacing $\tau$ by $\widetilde \tau$ has the effect of exchanging
$c$ with $c'$ and
$\nu$ with $\mu$.  In particular, since we have already shown that
$\sqrt{2+r}\in\calA_{-2,2}^{\varep\mu,\mu}$, it follows that
$\sqrt{2-r}=\sqrt{2+\widetilde r}\in\calA_{-2,2}^{\varep\nu,\nu}$.
\end{proof}

\subsection{Deterministic square roots.} \label{subsec:canonical}

Theorem~\ref{thm:main2} describes $\prod \calS_{k,\ell}^{\pm,\pm}$ in 
terms of square roots of $\tau$, $\tau+1$, or $\tau/(\tau+1)$.
In this section, we give a prescription for how this square root
is unambiguously determined. As a word of warning, these square roots do
not deserve to be called ``canonical''. For example, when $q=11$ and
$\tau=1$, our deterministic square root of $\tau$ turns out to be $-1$.

As usual, we set $j=4\tau/(\tau+1)$, $\ell=4/(\tau+1)$, $r=\ell-2$.
(See equations~(\ref{eq_tau})--(\ref{eq_tau2}).)
Assume $\tau \in\F_q \setminus \{0,-1\}$; then
$j$ and $\ell$ are nonzero, and $r\not \in \{2,-2\}$.
Write $r=\ang u$; then $u=(r\pm\sqrt{r^2-4\,\,})/2$ may be explicitly computed in a quadratic extension.  
Since $r\not \in \{2,-2\}$, we know $u^2\ne1$.
Any rational function of $u$ that is invariant under $u\mapsto 1/u$
is a symmetric function of $u$ and $1/u$, so may be expressed as
a rational function of $u+1/u$ and $u\cdot 1/u=1$. In other words,
it is a rational function of $u+1/u=r$.
Since $r=2(1-\tau)/(\tau+1)$, it is also a rational function of
$\tau$.  
Our goal is to produce well-determined
square roots of $\tau$, $\tau+1$, or $\tau/(\tau+1)$
when exactly one of these quantities is a square in $\F_q$.
We do this by finding a rational expression in $u$ whose square is $\tau$,
$\tau+1$, or $\tau/(\tau+1)$, and that is invariant under $u\mapsto 1/u$.

\noindent({\it i}) Suppose that $\tau \in \calA_{0,1}^{+-}$
(equivalently, $j$ and $\ell$ are nonsquares).
We will produce a deterministic square root of $\tau=j/\ell$.
By Lemma~\ref{lem:tauvw} or
Table~\ref{table:rho_ru}, $r \in \calA_{-2,2}^{-\varepsilon,-}$ and
$u^{(q-\varepsilon)/2}=-1$. Let
$$a_1 = \frac{u^{(q-\varepsilon)/4} - u^{-(q-\varepsilon)/4}}{u-1/u}.$$
This is invariant under $u \mapsto 1/u$, and 
\begin{equation*}
a_1^2 = \frac{u^{(q-\varepsilon)/2} + u^{-(q-\varepsilon)/2}-2}{u^2+1/u^2-2}
= \frac{(-1) + (-1) -2}{(u+1/u)^2-4} = \frac{-4}{r^2-4}.
\end{equation*}
Thus, $(2/a_1)^2 = 4-r^2=j\ell$, $(2/a_1\ell)^2 = j/\ell = \tau$.
We define our square root by $\sqrt{\tau\,} =\sqrt{j/\ell\,}= 2/(a_1\ell)$.

\rule[0mm]{0mm}{8mm} 
\noindent({\it ii})  Suppose that $\tau \in \calA_{0,1}^{-+}$
(equivalently, $j$ is a nonsquare and $\ell$ is a nonzero square). 
We will produce deterministic square roots of $\tau+1$ and $\ell$.
By Lemma~\ref{lem:tauvw} or
Table~\ref{table:rho_ru}, $r \in \calA_{-2,2}^{-\varepsilon,+}$ and
$u^{(q+\varepsilon)/2}=1$. Let $a_2=\ang{u^{(q-\varepsilon)/4}}$.
Then $a_2^2=u^{(q-\varepsilon)/2}+u^{-(q-\varepsilon)/2}+2$.
Now $u^{(q-\varepsilon)/2}=u^{(q+\varepsilon)/2-\varepsilon}=u^{-\varepsilon}$,
so $a_2^2=u^{-\varepsilon}+u^\varepsilon+2=r+2=\ell$. We define 
$\sqrt{\ell\,}=a_2$.
Since $\ell=4/(\tau+1)$, we define $\sqrt{\tau+1}=2/a_2=2/\sqrt{\ell\,}$.

\rule[0mm]{0mm}{8mm} 
\noindent({\it iii})  Suppose that $\tau\in\calA_{0,1}^{--}$
(equivalently, $j$ is a nonzero square and $\ell$ is a nonsquare).
We will produce deterministic square roots of $\tau/(\tau+1)$ and $j$.
By Lemma~\ref{lem:tauvw} or Table~\ref{table:rho_ru}, $r \in \calA_{-2,2}^{\varepsilon,-}$ and
$u^{(q+\varepsilon)/2}=-1$. Let $a_3=\ang{(-u)^{(q-\varepsilon)/4}}=\jacobi2q \ang{u^{(q-\varepsilon)/4}}$.
Then $a_3^2=\ang{u^{(q-\varepsilon)/4}}^2 = u^{(q-\varepsilon)/2}+u^{-(q-\varepsilon)/2}+2$.
Now $u^{(q-\varepsilon)/2}=u^{(q+\varepsilon)/2-\varepsilon}=-u^{-\varepsilon}$,
so $a_3^2=-u^{-\varepsilon}-u^\varepsilon+2=-r+2=j$. 
We define $\sqrt{j\,}=a_3$ to be our deterministic square root.
Since $j=4\tau/(\tau+1)$, we define $\sqrt{\tau/(\tau+1)}=a_3/2$.

\medskip
{\bf Remark.}\   Exchanging $j$ and $\ell$ is equivalent to $\tau \mapsto 1/\tau$, $r \mapsto -r$,
$u \mapsto -u$. Further, $\tau \in \calA_{0,1}^{\varepsilon_1,\varepsilon_2} \iff 1/\tau \in \calA_{0,1}^{\varepsilon_1,
\varepsilon_1 \varepsilon_2}$. Denote $\widetilde \tau = 1/\tau$, $\widetilde r = -r$, $\widetilde u = - u$.
In ({\it i}), if we set $\widetilde a_1 = (\widetilde u^{m}+\widetilde u^{-m})/
(\widetilde u - 1/\widetilde u)$ with $m=(q-\varepsilon)/4$, then 
$\widetilde a_1 = (-1)^{m+1} a_1 = -\jacobi 2q a_1$. Thus,
$\sqrt{\widetilde\tau} \sqrt\tau =2/(\widetilde a_1 j) \x 2/(a_1 \ell) = -\jacobi2q (2/a_1)^2/(j\ell) = - \jacobi2q$.
In ({\it ii}), $\tau \in \calA_{0,1}^{-+}$ implies $\widetilde \tau \in \calA_{0,1}^{--}$. 
The square root of $\ell$ is defined with respect to~$\tau$ in ({\it ii}) as $\sqrt {\ell\,}=a_2=\ang{u^m}$.
On the other hand, $\sqrt{\ell\,}$ is defined with respect to~$\widetilde \tau$ in ({\it iii}) as $\sqrt{\ell\,} = \widetilde a_3$,
where $\widetilde a_3 = \ang{\left(-\widetilde u\right)^m}=\ang{u^m} = a_2$. These two definitions for $\sqrt{\ell\,}$ are consistent.

\subsection{Proof of Theorem~\ref{thm:main1}({\it iv}) and Theorem~\ref{thm:main2}.} 
\label{subsec:main2Proof}
\parbox[c][20pt]{30pt}{ $\phantom{hello world}$ }

We are in position to prove Theorems~\ref{thm:main1}({\it iv}) and~\ref{thm:main2}.
As usual, let $\tau \in \F_q \setminus \{0,-1\}$,  and let $j,k,\ell,r,\tau'$ be as in (\ref{eq_tau})--(\ref{eq_tau2}).
Write $r=\ang {w^2}=\ang{u}$, where $w$ is as in Lemma~\ref{lem:tauvw} and $u=w^2$.
The hypothesis that $\tau \ne 0,\infty$ implies $u^2\ne 1$. 

\noindent{\it Proof of Theorem~\ref{thm:main1}({\it iv}):}\ 
The hypothesis is that $\tau\in\calA_{0,1}^{++}$.
By Lemma~\ref{lem:tauvw}, $r \in \calA_{-2,2}^{\varepsilon,+}$
and $u^{(q-\varepsilon)/2}=w^{q-\varepsilon}=1$.  
By (\ref{ST}) and Proposition~\ref{prop:Sklu}({\it i}),
$\prod\calT_{j,\ell}^{--} = \prod \calS_{k,\ell}^{-\varepsilon,-} =\prod \calS_{r-2,r+2}^{-\varepsilon,-}
=\jacobi2q \ang{u^{(q-\varepsilon)/4}}$.
By Proposition~\ref{prop:ru2} and (\ref{eq_ell}),  
$u^{(q-\varepsilon)/4}=\jacobi{2\pm2/\sqrt{\tau+1}}q = \jacobi{2\pm\sqrt{\ell\,}}q$. 
Thus,
$$\prod \calT_{j,\ell}^{--} 
=\jacobi2q \ang{u^{(q-\varepsilon)/4}}=\jacobi2q \jacobi{2\pm\sqrt{\ell\,}}q \cdot 2 = \jacobi{2\pm\sqrt{j\,}}q \cdot2,$$
where in the last step we use (\ref{jell}). 
By Proposition~\ref{prop:Tellj},
$\prod \calT_{\ell,j}^{--}= \jacobi2q \prod\calT_{j,\ell}^{--}$.
\qed

\noindent{\it Proof of Theorem~\ref{thm:main2}:}\  We continue with the
above notation. \\
({\it i}) The hypothesis is that $\tau\in\calA_{0,1}^{+-}$.
By Lemma~\ref{lem:tauvw} or Table~\ref{table:rho_ru}, $r\in \calA_{-2,2}^{-\varepsilon,-}$.
Then, by Proposition~\ref{prop:Sklu}({\it ii}),
\begin{equation*} \prod \calS_{k,\ell}^{\varepsilon,+}=\prod \calS_{r-2,r+2}^{\varepsilon,+} = -\jacobi 2q
\left(u^{(q-\varepsilon)/4}-u^{-(q-\varepsilon)/4}\right)/(u-1/u)
= -\jacobi2q a_1,
\end{equation*}
where $a_1$ is as in Section~\ref{subsec:canonical}. The deterministic square
root of $\tau$ is defined by $\sqrt\tau=2/(a_1\ell)$.  Thus,
$a_1=2/(\sqrt\tau \ell)=2/(\ell\sqrt{j/\ell\,\,})$, and the result follows.

\noindent ({\it ii})\ 
The hypothesis is that $\tau\in\calA_{0,1}^{-+}$. 
By Lemma~\ref{lem:tauvw} or Table~\ref{table:rho_ru}, $r\in \calA_{-2,2}^{-\varepsilon,+}$.
By Proposition~\ref{prop:Sklu}({\it i}), $\prod \calS_{r-2,r+2}^{-\varepsilon,-}=
\jacobi 2q \ang{u^{(q-\varepsilon)/4}}$. This is equal to $\jacobi2q a_2$,
where $a_2=\sqrt{\ell\,}=2/\sqrt{\tau+1}$ is as in Section~\ref{subsec:canonical}.
The result follows.

\noindent ({\it iii})\ 
We are assuming $\tau\in\calA_{0,1}^{--}$.
By Lemma~\ref{lem:tauvw} or Table~\ref{table:rho_ru}, $r\in \calA_{-2,2}^{\varepsilon,-}$ 
and $u^{(q+\varepsilon)/2}=-1$.
By Proposition~\ref{prop:Sklu}({\it i}), $\prod \calS_{k,\ell}^{-\varepsilon,-}=\prod \calS_{r-2,r+2}^{-\varepsilon,-}=
\jacobi2q \ang{u^{(q-\varepsilon)/4}}$. This equals $a_3=\sqrt{j\,\,}= 2\sqrt{\tau/(\tau+1)}$,
where $a_3$ is defined in Section~\ref{subsec:canonical}.  
\qed

\section{Rational reciprocity} \label{sec:reciprocity}

It sometimes happens that more than one row of Table~\ref{tableTjl1}
applies for a given value of $\tau$. It turns out that one can obtain
some known results about rational reciprocity in this way.
An excellent reference for reciprocity laws, including their
historical development, is by Lemmermeyer, \cite{Lemmermeyer}.

For example, if $\tau=1$ 
then $\tau$ and $\tau+1$ are squares iff $\jacobi 2q=1$,
or equivalently, iff $q\equiv \pm 1 \pmod 8$. Assuming this, then the
third row of Table~\ref{tableTjl1} tells us that 
\begin{equation} \prod \calT_{2,2}^{-+} = 
(-1)^{\lfloor (q+3)/8 \rfloor}. \label{comp1}
\end{equation}
On the other hand,
the two bottom rows of Table~\ref{tableTjl1} tell us that
$\prod \calT_{2,2}^{-+}= \jacobi{1\pm 1/\sqrt2}q$.
On multiplying the right side by $\jacobi 2q=1$, we obtain
\begin{equation} \prod \calT_{2,2}^{-+}= \jacobi{2\pm \sqrt2}q. \label{comp2} 
\end{equation}
On comparison of (\ref{comp1}) and~(\ref{comp2}), we find 
\begin{equation} \jacobi{2 \pm \sqrt 2}q = \begin{cases} 
(-1)^{(q-1)/8} & \text{ if $q\equiv 1 \pmod 8$,} \\
(-1)^{(q+1)/8} & \text{ if $q\equiv 7 \pmod 8$.} 
\end{cases} \label{biquad2}
\end{equation}
A different proof of this result can be found in \cite[p.~166]{Lemmermeyer}.

Another interesting phenomenon is the following theorem and corollary.

\begin{proposition} \label{prop:beyond}
Suppose that $\jacobi2q=1$ (equivalently, $q\equiv \varepsilon \pmod 8$ where $\varepsilon \in \{1,-1\}$),
and let $\sqrt2\in\F_q$ denote a (fixed) square root of 2. If $q \equiv \varepsilon \pmod{16}$ then
$$\prod \calT_{2-\sqrt2,2+\sqrt2}^{--} = 
\prod \left\{ a \in \F_q^\x : \text{$2\pm(\sqrt2+ a)$ are nonsquares} \right\}
= (-1)^{(q-\varepsilon)/16}\cdot 2.$$
If $q \equiv 8 + \varepsilon \pmod{16}$, then
$$\prod \calT_{2-\sqrt2,2+\sqrt2}^{++} = \prod \left\{ a \in \F_q^\x : 
\text{$2\pm(\sqrt2+ a)$ are nonzero squares}
\right\}
= (-1)^{(q+8-\varepsilon)/16}\cdot \sqrt2.$$
\end{proposition}

\begin{proof} 
Let $\zeta$ be a primitive eighth root of unity in $\cj\F_q$, and let $r =\ang{ \zeta}$.   
Now $\zeta^{q-\varepsilon}=1$,
so $r^q = \ang{\zeta^q}=\ang{\zeta^{\varepsilon}} = \ang\zeta = r$.  Further, since
$\zeta^2=i$ is a square root of $-1$, we have $r^2=\zeta^2+\zeta^{-2}+2
=i+1/i+2=2$; therefore $r=\pm \sqrt2$.  
By changing $\zeta$ to $-\zeta$ if necessary, we can assume 
$r=\sqrt2=\ang{\zeta}$. 

If $q\equiv \varepsilon \pmod{16}$ then by (\ref{biquad2}), $2\pm \sqrt2$ 
are squares, so
$\sqrt2 \not \in \calA_{-2,2}^{-\varepsilon,-}$. By Proposition~\ref{prop:Sklu}({\it i}) and (\ref{ST}),
$\prod \calT_{2-\sqrt2,2+\sqrt2}^{--}= \ang{(-\zeta)^{(q-\varepsilon)/4}}$.
Writing $q=\varepsilon+16\mu$ where $\mu\in\Z$, the right side is 
$\ang{(-\zeta)^{4\mu}}= (-1)^\mu \cdot 2$.  This proves the first formula.

Next, if $q\equiv 8 + \varepsilon \pmod{16}$ then by (\ref{biquad2}),
$2\pm \sqrt2$ are nonsquares, so 
$\sqrt2 \not \in \calA_{-2,2}^{\varepsilon,+}$. By Proposition~\ref{prop:Sklu}({\it ii}) and (\ref{ST}),
$$\prod \calT_{2-\sqrt2,2+\sqrt2}^{++}
=-\left((-\zeta)^{(q-\varepsilon)/4}-(-\zeta)^{-(q-\varepsilon)/4}\right)/(\zeta-1/\zeta).$$
Writing $q = 8 + \varepsilon+16\mu$, the numerator is 
$$(-1)((-\zeta)^{2+4\mu}-(-\zeta)^{-(2+4\mu)}).$$
Since $\zeta^4=-1$, this equals
$(-1)^{\mu+1} (\zeta^{2}-\zeta^{-2})$.   Thus,
$\prod \calT_{2-\sqrt2,2+\sqrt2}^{++}=(-1)^{\mu+1}(\zeta^2-\zeta^{-2})/(\zeta-1/\zeta)=(-1)^{\mu+1}\ang{\zeta}
=(-1)^{\mu+1} \sqrt2 = (-1)^{(q+8-\varepsilon)/16} \sqrt2$.
\end{proof}

\begin{corollary} \label{cor:beyond} Suppose that $q \equiv \varepsilon \pmod{16}$, so that $2+\sqrt 2$ is a square in
$\F_q$ (for either choice of square root) by~(\ref{biquad2}). Then $2+\sqrt{2+\sqrt 2}$ and $2-\sqrt{2+\sqrt2}$ are both
squares if $q \equiv \varepsilon \pmod{32}$, and otherwise they are both nonsquares.
\end{corollary}

\begin{proof} They are both squares or both nonsquares because the product, $4-(2+\sqrt 2)=2-\sqrt 2$,
is a square.
In Table~\ref{tableTjl1}, set $j=2-\sqrt 2$ and $\ell=2+\sqrt 2$. 
Then $\tau = (2-\sqrt 2)/(2+\sqrt 2)$, $\tau+1 = 4/(2+\sqrt 2)$, $2+2/\sqrt{\tau+1}=2+\sqrt{2+\sqrt 2}$.
According to the table, $\prod \calT_{j,\ell}^{--}=2$ if $2+2/\sqrt{\tau+1}$ is a square, and $-2$ otherwise.
On the other hand,  Proposition~\ref{prop:beyond} asserts that 
$\prod \calT_{j,\ell}^{--}=(-1)^{(q-\varepsilon)/16}\cdot 2$.
The result follows.
\end{proof}

One can continue down this path indefinitely, although the easiest proof does not involve Dickson polynomials.
The next lemma is well known.

\begin{lemma} \label{lem:angRat} If $u \in \cj\F_q $ is a primitive $d$-th root of unity, where $(d,q)=1$,
then $\ang u \in \F_q \iff q \equiv \pm1 \pmod d$.
\end{lemma}

\begin{proof}  In general, $\ang u = \ang v$ iff $v \in \{u,1/u\}$. Thus,
$\ang u \in \F_q \iff \ang{u^q}=\ang u \iff u^q \in \{u,1/u\}$.  The latter
is equivalent to $u^{q-1}=1$ or $u^{q+1}=1$, \ie, $d|q-1$ or $d|q+1$.
Thus, $\ang u \in \F_q \iff q \equiv \pm1 \pmod d$.
\end{proof}

\begin{proposition} In $\cj\F_q$, let $b_0=\sqrt 2$, $b_1=\sqrt{2+b_0}$, $b_2=\sqrt{2+b_1}$, and in general,
$b_i = \sqrt{2+b_{i-1}}$, where the choice of square root is arbitrary 
at each stage.
Then $b_i \in \F_q$ iff $q \equiv \varepsilon \pmod{2^{i+3}}$, 
where $\varepsilon=\jacobi{-1}q \in \{1,-1\}$.
\end{proposition}
\begin{proof} First, we show by induction that $b_j = \ang{\zeta_{j+3}}$, where $\zeta_j$ is a 
primitive $2^j$ root of unity\footnote{Here the $2^{j+1}$ choices of square 
roots in the definition of $b_j$ correspond to the $2^{j+2}$ distinct choices 
for primitive $2^{j+3}$ roots of unity, modulo
the relation $\ang \zeta=\ang {\zeta^{-1}}$.}  
in $\F_q$.  
To start the induction, $\sqrt2=\ang{\zeta_3}$ was demonstrated
in the proof of Proposition~\ref{prop:beyond}.
Now assume the induction hypothesis for $j$ and we prove it for $j+1$.
Let $\zeta_{j+4}^2=\zeta_{j+3}$; then $\zeta_{j+4}$ is a primitive 
$2^{j+4}$ root of unity. 
Now 
$\ang{\zeta_{j+4}}^2 = \ang{\zeta_{j+3}}+2=b_j+2$, so $\ang{\zeta_{j+4}}
= \pm \sqrt{b_j+2} = \pm b_{j+1}$. 
Replacing $\zeta_{j+4}$ by its negative if needed, we may assume that
$\ang{\zeta_{j+4}}=b_{j+1}$.  Finally, $b_j\in\F_q \iff \ang{\zeta_{j+3}}\in\F_q
\iff q \equiv \pm 1 \pmod{2^{j+3}}$ by Lemma~\ref{lem:angRat}.
\end{proof}

The above results generalize in various ways.  

\begin{proposition}  \label{prop:b0}
In $\cj\F_q$, suppose $u_0$ is a primitive $2k$-th root of unity and 
$b_0 = \ang{u_0}$. Let $b_1=\sqrt{2+b_0}, \ldots, b_i = \sqrt{2+b_{i-1}},\ldots$, 
for any choice of square roots. Then \\
({\it i})\ There are $u_1,u_2,\ldots \in \cj\F_q$ such that $u_i^2=u_{i-1}$ and $b_i= \ang{u_i}$. \\
({\it ii})\ For $i\ge0$, $b_i \in \F_q$ iff $q \equiv \pm 1 \pmod{2^{i+1} k }$. \\
({\it iii})\ Suppose $q \equiv \varepsilon \pmod{4k}$,
where as usual $\varepsilon = \jacobi{-1}q$, and write $q=\varepsilon + 4k\mu$.  Then
$$\prod \calT_{2-b_0,2+b_0}^{--} = \prod\{a \in \F_q^\x : \text{$2\pm(b_0+a)$ are nonsquares} \} = (-1)^{(k+1)\mu}\cdot 2.$$
\end{proposition}
\begin{proof}  
Let $u_1^2=u_0$. Then
$\ang{u_1}^2=u_0+u_0^{-1} + 2 = b_0+2$, so $b_1 = \pm \ang{u_1}$.
Replacing $u_1$ by $-u_1$ if necessary, we can assume $b_1 = \ang{u_1}$.
By induction, we can find $u_2, u_3, \ldots$ such that $u_i^2=u_{i-1}$ and
$b_i = \ang{u_i}$. Clearly $u_i$ is primitive of order $2^i(2k)$.
By Lemma~\ref{lem:angRat}, $b_i \in \F_q \iff q \equiv \pm 1 \pmod{2^i(2k)}$. 
This proves ({\it i}) and ({\it ii}).
Now we prove ({\it iii}).  The hypothesis $q \equiv \varepsilon \pmod {4k}$
implies $b_1\in\F_q$ by ({\it ii}), therefore $2+b_0$ is a square.
In particular, $b_0 \not \in \calA_{-2,2}^{-\varepsilon,-}$. 
By Proposition~\ref{prop:Sklu}({\it i}) and (\ref{ST}), 
$\prod \calT_{2-b_0,2+b_0}^{--} = \ang{(-u_0)^{k\mu} }$.
Since $u_0^k=-1$, we have $(-u_0)^{k\mu} = (-1)^{k\mu}(-1)^\mu=(-1)^{(k+1)\mu}$.
The result follows.
\end{proof}

\begin{proposition} \label{prop:threeTree} ({\it i}) Suppose $(q,6)=1$.
In $\cj \F_q$, 
let $b_0=\sqrt{3}$, $b_1 = \sqrt{2+b_0}, \ldots,
b_i = \sqrt{2+b_{i-1}},\ldots$, for any choice of square roots.
Then $b_i \in \F_q$ iff $q \equiv \pm1 \pmod{12\cdot 2^{i}}$. \\
({\it ii}) Suppose $(q,10)=1$. In $\cj\F_q$, let $b_0 = (1-\sqrt5)/2$, $b_1=\sqrt{2+b_0},\ldots,
b_i = \sqrt{2+b_{i-1}},\ldots$, for any choice of square roots. Then 
$b_i \in \F_q$ iff $q \equiv \pm1
\pmod{10\cdot 2^i}$. 
\end{proposition}
\begin{proof}
Using Proposition~\ref{prop:b0}, it suffices to show that 
\begin{equation*}\text{$\sqrt 3 = \ang{\zeta_{12}}$ 
and $(1-\sqrt5)/2 = \ang{\zeta_{10}}$,} \end{equation*}
where $\zeta_j$ denotes a primitive $j$-th root of unity.

To prove this for $\sqrt3$, let $\omega
=-\zeta_{12}^2$; then $\omega^3=-\zeta_{12}^6=1$. Noting that 
$\omega^2+\omega+1=0$, we have $\ang{\omega}=-1$.
Then $\ang{\zeta_{12}}^2 = \ang{\zeta_{12}^2}+2=\ang{-\omega}+2 =-\ang{\omega}+2=3$.
Thus, $\sqrt3 = \pm \ang{\zeta_{12}} = \ang{\pm \zeta_{12}}$.
Replacing $\zeta_{12}$ by $-\zeta_{12}$ if necessary, we can assume $\sqrt3 = \ang{\zeta_{12}}$.

To prove this for $(1-\sqrt5)/2$, let $\zeta$ be a primitive fifth root of unity.
Now $\ang{\zeta} + \ang{\zeta^2}
=\zeta+\zeta^4+\zeta^2+\zeta^3=-1$, so $(\ang\zeta+1/2)^2 = \ang{\zeta}^2
+\ang\zeta + 1/4 = (\ang{\zeta^2}+2) + \ang\zeta+1/4=2-1+1/4=5/4$.
This shows that $2\ang{\zeta}+1$ is a square
root of 5. The other square root of 5 is $2\ang{\zeta^2}+1= -(2\ang\zeta+1)$.
Replacing $\zeta$ by $\zeta^2$ if necessary gives
$\ang{-\zeta} = -\ang\zeta = (1-\sqrt 5)/2$. Set $\zeta_{10}=-\zeta$, which is a primitive tenth root of unity.
\end{proof}

We remark that Lemma~\ref{lem:angRat}, together with the formulas 
\begin{equation}
\text{$\ang{\zeta_8}=\sqrt2$,\quad $\ang{\zeta_{12}}=\sqrt3$,\quad and\quad $\ang{\zeta_{10}}=(1-\sqrt 5)/2$,} 
\label{specialAngs}
\end{equation}
give one-line proofs: 
\begin{equation} \jacobi 2q = 1 \iff \sqrt2 \in\F_q \iff \ang{\zeta_8} \in \F_q \iff q\equiv \pm 1 \pmod 8,\label{oneline2}
\end{equation}
\begin{equation} \jacobi 3q = 1 \iff \sqrt3 \in\F_q^\x \iff \text{$3\nmid q$ and $\ang{\zeta_{12}} \in \F_q$} \iff q\equiv \pm 1 \pmod {12}.\label{oneline3}
\end{equation}
\begin{equation} \jacobi 5q = 1 \iff \sqrt5 \in\F_q^\x \iff \text{$5\nmid q$ and $\ang{\zeta_{10}} \in \F_q$} \iff q\equiv \pm 1 \pmod {10}.\label{oneline5}
\end{equation}

\begin{proposition} \label{prop:sqrt3} Suppose that $\jacobi3q=1$ (equivalently, $q\equiv \pm1 \pmod{12}$),
and let $\sqrt3 \in \F_q$ denote a (fixed) square root of~3. Let $\nu=(-1)^{(q-\varep)/12}$, \ie, $\nu=1$ if
$q \equiv \pm1 \pmod{24}$, $\nu=-1$ otherwise. Then
$$\prod \calT^{-\nu,-\nu}_{2-\sqrt3,2+\sqrt3} = (-1)^{\mu} \cdot 2,\qquad\text{where $\mu=\lfloor (q+1)/24 \rfloor$.}$$
\end{proposition}

\begin{proof} 
Write $\sqrt 3 = \ang{u}$, where $u=\zeta_{12}$ is a primitive twelfth root of unity.
By Proposition~\ref{prop:b0}({\it iii}), if $q\equiv \pm1 \pmod{24}$ then $\prod \calT_{2-\sqrt3,2+\sqrt3}^{--} = (-1)^{\mu}\cdot 2$.

Now suppose $q \not \equiv \pm 1 \pmod{24}$. Then we may write $q=\varepsilon+12+24\mu$, where $\varepsilon \in \{1,-1\}$.
By Proposition~\ref{prop:b0}({\it ii}), $\sqrt{2+\sqrt 3} \not \in \F_q$, \ie, $\jacobi{2+\sqrt3}q = -1$. Then
$\sqrt3 \not \in \calA_{-2,2}^{\varepsilon,+}$.
By Proposition~\ref{prop:Sklu}({\it ii}) and (\ref{ST}), 
$$\prod \calT_{2-\sqrt3,2+\sqrt3}^{++}= 
-\left((-u)^{(q-\varepsilon)/4}-(-u)^{-(q-\varepsilon)/4}\right)/(u-1/u).$$
Since $u^6=-1$, the numerator of the right side is
$(-1)((-u)^{3+6\mu}-(-u)^{-(3+6\mu)})=u^{3+6\mu}-u^{-3-6\mu}= (-1)^\mu(u^3-u^{-3})$. Then the right side is
$(-1)^\mu(u^3-u^{-3})/(u-1/u) = (-1)^\mu(u^{2}+1+u^{-2})$.
Let $\omega = -u^2$; then $\omega^3=-u^6=1$, so $\omega^2+\omega+1=0$. We have $u^2+1+u^{-2}=\ang{-\omega}+1=2$, 
so $\prod\calT_{2-\sqrt3,2+\sqrt3}^{++} =(-1)^\mu \cdot 2$.
\end{proof}

\begin{proposition} \label{prop:sqrt5} Suppose that $\jacobi5q=1$ (equivalently, $q\equiv \pm1 \pmod{5}$),
and let $r = (1-\sqrt5)/2$, where $\sqrt5 \in \F_q$ denotes a (fixed) square root of~5. 
Then
\begin{eqnarray*} \prod \calT_{2-r,2+r}^{--} &=& 2\qquad  \text{\quad\, if $q \equiv \pm1 \pmod{20}$,}  \\
\prod \calT_{2-r,2+r}^{++} &=& - \varepsilon r \qquad \text{otherwise.}
\end{eqnarray*}
\end{proposition}

\begin{proof} Write $r= \ang{u}$, where $u=\zeta_{10}$ is a primitive tenth root of unity.  
If $q \equiv \pm 1 \pmod{20}$  then applying
Proposition~\ref{prop:b0}({\it iii}) with $k=5$ gives $\prod \calT_{2-r,2+r}^{--} = 2$. 
Now suppose that $q \not \equiv \pm 1 \pmod{20}$.
Being that $q$ is odd and $q\equiv \pm1 \pmod 5$ by hypothesis, we have $q\equiv \pm1 \pmod{10}$. 
Write $q = 10-\varepsilon+20\mu$, where $\varepsilon\in \{1,-1\}$.
Note that $q\equiv \varepsilon \pmod 4$, so $\varepsilon = \jacobi{-1}q$.
By Proposition~\ref{prop:b0}({\it ii}) with $i=1$, $2+r$ is a nonsquare, so
$r \not \in \calA_{-2,2}^{\varepsilon,+}$. By Proposition~\ref{prop:Sklu}({\it ii}),
$$\prod \calT_{2-r,2+r}^{++} = - \left((-u)^{(q-\varepsilon)/4}-(-u)^{-(q-\varepsilon)/4}\right)/(u-1/u).$$
We have $(q-\varepsilon)/4=(10-2\varepsilon+20\mu)/4=2+5\mu + (1-\varepsilon)/2$. Let $\zeta=-u$, so $\zeta^5
=-u^5=1$. Then $(-u)^{(q-\varepsilon)/4} = \zeta^{2+5\mu+(1-\varepsilon)/2}$. This is $\zeta^2$ if $\varepsilon=1$,
or $\zeta^3=\zeta^{-2}$ if $\varepsilon =-1$. Therefore,
$$\prod \calT_{2-r,2+r}^{++}= \begin{cases} 
(\zeta^2-\zeta^{-2})/(\zeta-\zeta^{-1}) & \text{if $\varepsilon=1$} \\
(\zeta^{-2}-\zeta^{2})/(\zeta-\zeta^{-1}) & \text{if $\varepsilon=-1$}. 
\end{cases} $$
Thus, $\prod \calT_{2-r,2+r}^{++}= \varepsilon\ang\zeta = -\varepsilon \ang u = -\varepsilon r$.
\end{proof}

\end{document}